\documentclass[10pt]{article}
\usepackage{amsmath,amsthm,amsfonts,amssymb,bbm,
latexsym,cite}
\usepackage{psfrag,graphicx, appendix}
\numberwithin{equation}{section}
\newtheorem{theo}{Theorem}[section]
\newtheorem{prop}{Proposition}[section]
\newtheorem{lemme}{Lemma}[section]
\newtheorem{conjecture}{Conjecture}[section]
\newtheorem{definition}{Definition}[section]
\newtheorem{fact}{Fact}[section]
\newtheorem{claim}{Claim}[section]

\theoremstyle{remark}

\newtheorem{remark}{Remark}[section]

\newcommand{\bl}{\begin{lemme}}
\newcommand{\el}{\end{lemme}}

\newcommand{\brem}{\begin{remark}}
\newcommand{\erem}{\end{remark}}
\newcommand{\bconj}{\begin{conjecture}}
\newcommand{\econj}{\end{conjecture}}
\newcommand{\bdefi}{\begin{definition}}
\newcommand{\edefi}{\end{definition}}
\newcommand{\bt}{\begin{theo}}
\newcommand{\bfa}{\begin{fact}}
\newcommand{\efa}{\end{fact}}

\newcommand{\et}{\end{theo}}
\newcommand{\bp}{\begin{prop}}
\newcommand{\ep}{\end{prop}}
\newcommand{\be}{\begin{equation}}
\newcommand{\mbE}{\mathbb{E}}
\newcommand{\tTr}{{\rm{Tr}}}
\newcommand{\ee}{\end{equation}}
\newcommand{\mbC}{\mathbb{C}}
\newcommand{\mbP}{\mathbb{P}}
\newcommand{\mbR}{\mathbb{R}}

\renewcommand{\Re}{\mathrm{Re}}
\renewcommand{\Im}{\mathrm{Im}}
\title{Universality in the bulk of the spectrum for complex sample covariance matrices}
\author{ S. P\'ech\'e\thanks{Institut Fourier, 100 Rue des Maths,
    BP 74, F-38402 St Martin d'Heres. E-mail: Sandrine.Peche@ujf-grenoble.fr}}
\begin{document}
\maketitle
\begin{abstract}
We consider complex sample covariance matrices $M_N=\frac{1}{N}YY^*$ where $Y$ is a $N \times p$ random matrix with i.i.d. entries $Y_{ij}, 1\leq i\leq N, 1\leq j \leq p$ with distribution $F$. Under some regularity and decay assumption on $F$, we prove universality of some local eigenvalue statistics in the bulk of the spectrum in the limit where 
$N\to \infty$ and $\lim_{N \to \infty}p/N =\gamma$ for any real number $\gamma \in (0, \infty)$.  
\end{abstract}

\section{Introduction}
\subsection{Model and result}
This paper is concerned with universal properties of large complex sample covariance matrices in the bulk of the spectrum.
We consider $N \times p$ random matrices $Y=(Y_{ij})$ where the $Y_{ij}$'s are i.i.d. random variables with some probability distribution $F$.
Let then $M_N$ be the sample covariance matrix:
$$M_N=\frac{1}{N}YY^*.$$
In the whole paper we assume that $p\geq N$ and that $$\exists \gamma\in [1, \infty)\text{ such that }p/N \to \gamma \text{ as }N \to \infty.$$ The case where $p<N$ can be deduced from the above setting using the fact that $YY^*$ and $Y^*Y$ have the same non zero eigenvalues.
In the sequel we call $\lambda_1\geq \lambda_2\geq \cdots \geq \lambda_N$ the ordered eigenvalues of $M_N$ and $\pi_N:=\frac{1}{N}\sum_{i=1}^N \delta_{\lambda_i}$ its spectral measure.
Assuming that $F$ has a finite variance $\sigma^2$, Mar\v{c}enko and Pastur (\cite{MP}, see also \cite{Sil}) have shown that $\pi_N$ almost surely converges as $N \to \infty$ to the so-called Mar\v{c}enko-Pastur distribution $\rho^{MP}_{\gamma,\sigma^2}.$ This probability distribution depends on $\sigma^2$ only and not on any other detail (higher moments e.g.) of $F$: in this sense it is \emph{universal}. It is defined by the density with respect to the Lebesgue measure:
\be \label{densMP}\frac{d\rho^{MP}_{\gamma,\sigma^2}(x)}{dx}=\frac{\sqrt{(u_+-x)(x-u_-)}}{2\pi x \sigma^2}1_{x \in [u_-, u_+]},\ee
where $u_+=\sigma^2(1+\sqrt{\gamma})^2$ and $u_-=\sigma^2(1-\sqrt \gamma)^2.$

It has been conjectured (see \cite{Me} for instance) that, in the large-$N$-limit, some finer properties of the spectrum are also universal. For instance, the spacing between nearest neighbor eigenvalues in the vicinity of a point $u\in (u_-, u_+)$ is expected to be universal, under the sole assumption that the variance of $F$ is finite.
The spacing is actually believed to be ``more universal'' than the limiting Mar\v{c}enko-Pastur distribution, in the sense that it is expected to be the same as for Hermitian Wigner matrices.  
To investigate such local properties of the spectrum, we introduce the so-called \emph{local eigenvalue statistics} in the bulk of the spectrum. 
Given a symmetric function $f \in L^{\infty}
(\mathbb{R}^m)$ ($m$ fixed) with compact support, a point $u \in
[u_-,u_+],$ and a scaling factor $\rho_N$, we define the local eigenvalue statistic $S_N^{(m)}(f,u, \rho_N)$ by
\begin{equation}
S_N^{(m)} (f,u, \rho_N)= \sum_{i_1, \ldots, i_m}f(\rho_N(\lambda_{i_1}-u),\ldots,\rho_N(\lambda_{i_m}-u)),\label{les}
\end{equation}
where the sum is over all distinct indices from $\{1,\ldots, N\}.$
When $u$ is in the bulk of the spectrum, that is $u\in (u_-, u_+)$, the natural choice for the scaling factor is $\rho_N= N \rho^{MP}_{\gamma, \sigma^2}(u)$,
which should give the scale of the spacing between nearest neighbor eigenvalues in the vicinity of $u$.

\paragraph{}
We here prove universality of some local linear statistics in the bulk of the spectrum for a wide class of complex sample covariance matrices. We follow the approach used by \cite{SchleinYau3} (and a series of papers \cite{SchleinYau0}, \cite{SchleinYau00}, \cite{SchleinYau1}) where universality in the bulk of Wigner matrices is proved.
We now define the class of matrices under consideration in this paper.
Let $\mu$ be the real Gaussian distribution with mean $0$ and variance $1/2$.
Let $F$ be a complex probability distribution whose real and imaginary parts are of the form
\be \label{12prime}\nu(dx) = e^{-V(x)}\mu(dx) ,\ee
for some real function $V$ satisfying the following assumptions:\\
- $V \in \mathcal{C}^6$ and there exists an integer $k\geq 1$ such that
\be \sum_{j=1}^6|V^{(j)}(x)| \leq C_*(1 + x^2)^k, \label{(1.5)}\ee
- there exist $\delta_1, C' >0$ such that $\forall x\in \mbR$, \be \label{(1.6)}\nu(x)\leq C'e^{-\delta_1|x|^2}.\ee
This assumption can actually be relaxed as explained in Remark \ref{rem: expdecay} below to consider distributions $\nu$ with sub-exponential decay only.
\\- $F$ is normalized so that \be \label{(1.8)}\int x d\nu(x)=0 \text{ and }\int |x^2|d\nu(x)=1/2.\ee
In the sequel we consider sample covariance matrices $M_N=\frac{1}{N}YY^*$ where $Y=(Y_{ij})$ is a $N \times p$ random matrix such that:
\be \label{assumptionY}Y_{ij}, 1\leq i\leq N; 1\leq j\leq p\text{ are i.i.d. random variables with distribution $F$}. 
\ee
One shall remark that the condition (\ref{(1.8)}) can always be achieved by rescaling and does not impact on the generality of our next results.\\
We now give our two main results.
Let $\epsilon >0$ small be given.
\bt
\label{theo: uni2}
Assume that $F$ satisfies (\ref{12prime}) to (\ref{(1.8)}).
Let also $u\in (u_-+\epsilon, u_+-\epsilon)$ be a point in the bulk of the spectrum and $\rho_N=N\rho^{MP}_{\gamma, 1}(u).$ 
Then $$\lim_{N \to \infty}\mbE S_N^{(2)} (f,u, \rho_N)=\int_{\mbR^2}f(x,y) \Big (1-\left (\frac{\sin (x-y)\pi}{\pi(x-y)}\right)^2\Big )dx dy.$$
\et
\brem It is also possible to prove universality of local eigenvalue statistics for $m=3$ using the approach developed hereafter. Nevertheless to consider higher integers $m$, one needs to increase the regularity of $V$ (see Remark 1.1 in \cite{SchleinYau3}). 
\erem
We can also prove that the spacing distribution close to a point $u $ in the interior of the support
of Mar\v{c}enko-Pastur's law is universal.
Let $s \geq 0$ and $(t_N)$ be a sequence such that 
$\lim_{N \rightarrow \infty}t_N=+\infty$ and $\lim_{N \rightarrow \infty}\frac{t_N}{N^{1-\beta}}=0 $ for some $\beta>0.$
Let $u\in (u_-+\epsilon, u_+-\epsilon)$ be given and $\rho_N=N\rho^{MP}_{\gamma, 1}(u).$ 
Define then the ``spacing function'' of the eigenvalues by
\begin{equation}\label{spacing}
S_N(s,u):=
\frac{1}{2t_N}\sharp \left \{ 1\leq j \leq N-1,\, \lambda_{j+1}-\lambda_j
\leq \frac{s}{\rho_N},\,| \lambda_j-u | \leq \frac{t_N}{\rho_N}  \right \}. 
\end{equation}
Intuitively the expectation of the spacing function is the probability, knowing that there exists
an eigenvalue in the interval $[u-t_N, u+t_N]$, to find its nearest neighbor
within a distance $\frac{s}{\rho_N}.$
Finally, we define
\begin{equation}
p(s)=
\frac{d^2}{ds^2}\det(I-K)_{L^2(0,s)}, \text{ where } K(x,y):=\frac{\sin \pi(x-y)}{\pi(x-y)}. \label{p(s)}
\end{equation}

\bt \label{theo: spacing} Assume that $F$ satisfies (\ref{12prime}) to (\ref{(1.8)}).
Let also $u\in (u_-+\epsilon, u_+-\epsilon)$ be a point in the bulk of the spectrum. 
Then,
\begin{equation}
\lim_{N\rightarrow \infty}\mathbb{E}S_N(s, u)=\int_0^s p(w)dw. \end{equation}
\et

\brem \label{rem: expdecay}Using truncation and centralization techniques, it is possible to prove both Theorem \ref{theo: uni2}
and Theorem \ref{theo: spacing} when Assumption \ref{(1.6)} is replaced by the weaker assumption
$$\exists \, C_1, C_2>0 \text{ such that } \nu(x) \leq C_1e^{-C_2|x|}.$$ This extension is examined in full detail in Section 5 of \cite{SchleinYau3} for Wigner random matrices and readily extends to sample covariance matrices.

\erem

The first proofs of universality in the bulk of the spectrum of large random matrices have been obtained for the so-called invariant ensembles (\cite{D}, \cite{DKMVZ1}, \cite{DKMVZ2}). Their proof relies on the fact that the joint eigenvalue density of such ensembles can be computed and the asymptotic local eigenvalue statistics can then be determined. A breakthrough in the proof of the conjecture was obtained in \cite{Johansson} (following the idea of \cite{BHik} and \cite{BHik2}), proving universality in the bulk for the so-called Dyson Brownian motion model (\cite{Dy}). This then allowed to extend universality results to a wide class of non invariant Hermitian random matrix ensembles- the so-called Gauss divisible ensembles (see Subsection \ref{subsec: id} for the definition). \cite{GBAPecheCPAM} have then obtained the same universality result for Gauss divisible complex sample covariance matrices.
Very recently, such universality results have been greatly improved by \cite{TV} and \cite{SchleinYau3} for Hermitian Wigner random matrices. Both the papers remove the Gauss divisible assumption and only assume sufficient decay of the entries of Wigner random matrices. The approach of \cite{TV} is to make a Taylor expansion of local eigenvalue statistics of Wigner matrices. The core of the proof is then to show that these statistics depend, in the large-$N$-limit, only on the first four moments of the entries of the Wigner matrix. Proving that any Wigner matrix can be matched to a Gauss divisible matrix with the same first $4$ moments allows to prove a very general universality result. On the other hand, the approach of \cite{SchleinYau3} is to show that any Wigner matrix (under suitable decay of the entries) can be sufficiently well approximated by a Gauss divisible random matrix, so that the bulk universality follows. We also refer the reader to \cite{SchleinTao} where the two approaches are combined. \\
While writing this paper, a proof of universality in the bulk of the spectrum of large sample covariance matrices has been obtained in \cite{TV_SC} in the sole case where $p-N=O(N^\frac{43}{48})$, but with much milder assumptions on the decay of the distribution $\nu$. Therein $\nu$ satifies $\int |x|^{C_o}d\nu(x)<\infty$ for some sufficiently large $C_o$ in place of our assumptions (\ref{(1.5)}) and (\ref{(1.6)}). Their approach is based on the ideas developed in \cite{TV}.\\
Our paper closely follows the ideas developed in \cite{SchleinYau3}. We give an overview of the proof in the next subsection.
\subsection{\label{subsec: id}The idea of the proof }
Following the pioneering work of \cite{Johansson} and \cite{BHik}, \cite{GBAPecheCPAM} have shown universality of local eigenvalue statistics in the bulk of the spectrum for complex sample covariance matrices when the distribution of the sample is \emph{Gauss divisible}. We recall that a complex distribution $\mu_G$ is Gauss divisible if there exist a complex probability distribution $P$ and a non trivial complex Gaussian distribution $G$ such that $\mu_G=P\star G.$
Equivalently \cite{GBAPecheCPAM} consider random matrices of the form 
\be \widetilde{M}_N=\frac{1}{ N}\left ( W+	aX\right)\left ( W+	aX\right)^*\label{defense},\ee
where $W$ and $X$ are independent $N \times p$ complex random matrices both with i.i.d. entries. The $W_{ij}$'s are $P$ distributed and the $X_{ij}$'s are complex standard Gaussian random variables. In the above context $a$ is real number independent of $N$. 
The proof of \cite{Johansson} and \cite{GBAPecheCPAM} relies on three main steps: \\
-conditionnally on $H=\frac{1}{N} WW^*$, the eigenvalue process induced by $\widetilde{M}_N$ defines a so-called determinantal random point field;\\
-the corresponding correlation kernel can be expressed as a double integral in the complex plane depending on $H$ through its sole spectral measure $\mu_N$;\\
-under suitable assumptions on $W$ and thanks to concentration results for the spectral measure $\mu_N$ established by \cite{BaiConvRates}, \cite{BaiImpCR} and \cite{GuiZei}, the asymptotic analysis of the correlation kernel (and local statistics) can be performed.\\
The parameter $a$ is to be seen as the ``order of the Gaussian regularization'' of $P$.
In principle this result would yield a full proof of the universality conjecture if one could let $a$ approach (and be smaller than) $1/\sqrt N$. Unfortunately this idea fails whatever sharp concentration results can be established for $\mu_N$. The asymptotic analysis is not stable in this scale.
\paragraph{}A breakthrough to overcome this difficulty is obtained in \cite{SchleinYau3}. In \cite{SchleinYau0}, \cite{SchleinYau00}, \cite{SchleinYau1} concentration results are deeply sharpened so as to be able to consider a Gaussian regularization of order $a >>\frac{1}{\sqrt N}$. 
Given an arbitrary (non Gauss divisible) distribution $F$,
the main point is however to show that one can find a Gauss divisible distribution approximating $F$ sufficiently well so that one can deduce universality in the bulk for a sample covariance matrix with i.i.d. $F$ distributed entries.  
This is the main step achieved in \cite{SchleinYau3}, which we now briefly expose.

Consider the Ornstein-Uhlenbeck (OU) process
$$L=\frac{1}{4}\frac{\partial^2}{\partial x^2}-\frac{x}{2}\frac{\partial}{\partial x}, \quad \partial_t u=L u.$$
Given a complex probability distribution $\widetilde{\mu}$,  
let then $X_t$ be the $N \times p$ matrix process with initial distribution $\widetilde{\mu}^{\otimes Np}$ and whose entries (real and imaginary parts) evolve independently according to the OU process. Then $X_t$ evolves as 
$$t\mapsto e^{-t/2}\left ( \hat H+ (e^t-1)^{1/2}X\right),$$
where $\hat H$ has i.i.d. entries with distribution $\widetilde{\mu}$ and $X$ has standard complex Gaussian entries.
Let $e^{t\mathcal{L}}:=(e^{tL})^{\otimes Np}$ denote the dynamics of the OU process for all the matrix elements. 
Here, for small $t$, one should think that $t\simeq a^2.$ \\
\cite{SchleinYau3} prove that there exists a Gauss divisible distribution $\mu^t_G$ approximating $F$ in the total variation norm in a sufficiently good way. Roughly speaking, let 
$$\mu^t_G= e^{tL}(1-tL+tL^2/2) F^c, \quad \widetilde{\mu}=(1-tL +t^2L^2/2)F^c$$ where $F^c$ is obtained from $F$ after truncation and centralization. 
Intuitively, it is reasonable to expect that $\mu^t_G$ is a good approximation of $F$ in the scale $t<<1.$
For sake of completeness we here recall their result (Proposition 2.1). 

Let $\theta$ be a smooth cutoff function satisfying $\theta(x) = 1$ if $|x| \leq 1$ and $\theta(x) = 0$ for $|x| \geq 2$.
\bp \label{prop: approxF}
Let $V$ satisfy (\ref{(1.5)}), (\ref{(1.6)}) for some $k$ and (\ref{(1.8)}). Let $\lambda>0$ be
sufficiently small and $t=N^{\lambda-1}$. Let $c_N, d_N$ be real numbers so that $vd\mu$ defines a centered probability density if the function $v$ is given by:
$$v(x) := e^{-V^c(x)}, \,  V^c(x) := V (x)\theta\left ((x-c_N)N^{-\lambda/4k}\right ) + d_N.$$
Let then $\mathcal L:= L^{\otimes Np}, f_{N,p}:=(e^{-V})^{\otimes Np},$ and $f_{c, N,p}:=v^{\otimes Np}.$
\begin{enumerate} \item There exist constants $C>0$, $c>0$ depending on $k$ and $\lambda$ such that $\int |f_{c,N,p} - f_{N,p}| d\mu^{\otimes Np}\leq  C e^{-c(Np)^{c/2}}.$
\item  $g_t := (1-tL+t^2L^2/2)v$ is a probability measure with respect to $d\mu$. Setting $G_t := [g_t]^{\otimes Np},$ there exists a constant $C$ depending on $\lambda$ and the constants $C_*$ and $\delta_1$ defined in (\ref{(1.5)}) and (\ref{(1.6)}) such that $$\int \frac{|e^{t\mathcal{L}}G_t-f_{c,N,p}|^2}{e^{t\mathcal{L}}G_t}d\mu^{\otimes Np}\leq CNp t^{6-\lambda}\leq CN^{-4+8\lambda}\frac{p}{N}.$$
\end{enumerate}
\ep
The idea is then to prove Theorem \ref{theo: uni2}
and Theorem \ref{theo: spacing} for the Gauss divisible ensembles with small parameter $a\sim N^{(\lambda-1)/2}$ for some $\lambda>0$. Then, using Proposition \ref{prop: approxF}, and following the idea of \cite{SchleinYau3} Section 4, one can extend universality of local eigenvalue statistics $S_N^{(m)} (f,u, \rho_N)$ with $m=2,3$ and that of the spacing function to sample covariance matrices satisfying (\ref{assumptionY}).

The paper is organized as follows. In Section \ref{Sec: GD} we study eigenvalue statitics in the bulk of the spectrum for Gauss divisible sample covariance matrices. To this aim, we first recall some properties of the \emph{Deformed Wishart ensemble}. This is the conditional distribution of $\widetilde{M}_N$ knowing $W$. Such an ensemble is in particular known to be determinantal, as we recall. We then establish some improved convergence rates for the spectral measure of sample covariance matrices whose entries have a sub-Gaussian tail. These concentration results then allow to compute the asymptotic correlation functions as $N \to \infty$ in the regime where $a\to 0, a>>\frac{1}{\sqrt N}.$ We then prove Theorem \ref{theo: spacing}
and Theorem \ref{theo: uni2} in Section \ref{Sec: pfth}.\\
In the whole paper, we use $C$ and $c$ to denote constants whose value may vary from line to line.
\section{The Gauss divisible ensemble \label{Sec: GD}}
In this section we establish some universality results for the following Gauss divisible ensemble:
let \be \widetilde{M}_N=\frac{1}{ N}\left ( W+a X\right)\left ( W+a X\right)^*\label{defensem},\ee
be a Hermitian $N\times N$ random matrix where

\begin{itemize}
\item{$(H_0)$} $X$ and $W$ are independent $N \times p$ random matrices where $p\geq N$ and there exists $\gamma \geq 1$ so that $\lim_{N \to \infty}p/N =\gamma.$ 
\item{$(H_1)$}  $a=a_N=\sqrt{\frac{N^{\lambda}}{N}}$ where $\lambda>0$ is a (small) real number;
\item{$(H_2)$}  $X$ is a $N\times p$ random matrix with complex standard Gaussian entries: 
$$\Re X_{ij}, \Im X_{ij} \sim \mathcal{N}(0, 1/2), \forall\, 1\leq i\leq N,\, 1\leq j\leq p;$$ 
\item{$(H_3)$} $W=\left( W_{ij}\right), 1\leq i\leq N,\, 1\leq j\leq p,$ is a complex random matrix such that $\Re W_{ij},\, \Im W_{ij},\, 1\leq i \leq N, 1\leq j \leq p$ are i.i.d. and satisfy :
$$({\bf A_1}) \text{ There exists a constant $\delta_o>0$ such that }\mathbb{E}e^{\delta_o |W_{ij}|^2}<\infty, \,\forall \, i,j.$$
\end{itemize}

Without loss of generality, we also make the assumption that 
\be \label{varW}\mbE |W_{ij}|^2=1/4.\ee
Note that (\ref{varW}) can always be achieved by rescaling $\widetilde{M}_N$. 
From now on, we denote by $y_1\geq y_2\geq \cdots \geq y_N$ the ordered eigenvalues of $H=H_N=WW^*/N$ and let $\mu_N=\frac1N \sum_{i=1}^N \delta_{y_i}$ be its spectral measure.
Under assumption (\ref{varW}), it is known that $\mu_N$ almost surely converges to the Mar\v{c}enko-Pastur distribution with parameters $\gamma$ and $ 1/4$, $\rho^{MP}_{\gamma, 1/4}=\rho^{MP}_{\gamma}$, whose density function is given by (\ref{densMP}) with $\sigma^2=1/4.$ When $\gamma=1$, we simply denote $\rho^{MP}_1$ by $\rho^{MP}$ for short. 

\subsection{Correlation functions }
The sample covariance Gauss divisible ensemble has a nice mathematical structure that we are going to make use of in the sequel: the conditional distribution of $\widetilde{M}_N$ with respect to $W$ is the so-called complex deformed Wishart ensemble.
Such an ensemble has been widely studied in random matrix theory as it induces a determinantal random point field.
We recall some results used in \cite{GBAPecheCPAM} (see references therein) that will be needed for the sequel.
Let then $P_N^H(\lambda_1, \lambda_2, \ldots, \lambda_N)$ be the joint eigenvalue distribution induced by the conditional distribution of $\widetilde{M}_N$ w.r.t. $H$.
Then (see Section 3 in \cite{GBAPecheCPAM}), $P_N^H$ is absolutely continuous with respect to Lebesgue measure on $\mbR_+^N$. Its density $f_N^H$ is given by:
\begin{equation}f_N^H(x_1,x_2,\ldots,x_N)=
\frac{V(x)}{V(y)} \det
 \Bigl (\frac{1}{S}
 e^{\{-\frac{y_i+x_j}{S}\}}\,I_{\nu}\left (\frac{2\sqrt{y_i x_j}}{S}\right )\:
  \left (\frac{x_j}{y_i}\right )^{\frac{\nu}{2}}\Bigr )_{i,j=1}^N,
\label{densdefwish} 
\end{equation}
where $V(x):=\prod_{1\leq i<j\leq N}(x_i-x_j)$, $\nu=p-N$, $I_{\nu}$ is the modified Bessel function of the first kind, and $S=a^2/N$. The main tools to study local eigenvalue statistics are the so-called eigenvalue correlation functions (see Section \ref{Sec: pfth}). They are defined for any integer $1\leq m\leq N$ by
$ R_N^{(m)}(x_1, x_2, \ldots, x_m; H)=\frac{N!}{(N-m)!}\int_{\mbR_+^{N-m}}f_N^H(x_1,x_2,\ldots, x_m, \lambda_{m+1}, \ldots, \lambda_N)\prod_{i=m+1}^Nd\lambda_i.$ 
Then, for any integer $1\leq m\leq N$, one has that $$R_N^{(m)}(x_1, x_2, \ldots, x_m;H)=\det \left ( K_N(x_i,x_j;H)\right)_{i,j=1}^m,$$
for some correlation kernel $K_N(u,v;H)$. This gives the determinantal random point field structure. Furthermore $K_N(u,v;H)$ is given by 
\begin{eqnarray}\label{defKN}
&K_N(u,v;H):=&\frac{1}{2i^2\pi^2 S^2}\int_{i\mbR+A}dw\int_{\Gamma}dz e^{\frac{-z^2+w^2}{S}}K_{\nu}\Big (\frac{2w\sqrt v}{S}\Big )I_{\nu}\Big (\frac{2z\sqrt u}{S}\Big )\crcr
&& \prod_{i=1}^N \frac{w^2-y_i}{z^2-y_i}\left (\frac{w}{z}\right)^{\nu} \left (\frac{w+z}{w-z}-\frac{w-z}{w+z}\right),\end{eqnarray}
where $K_{\nu}$ is the modified Bessel function of the second kind. $\Gamma$ is a contour oriented counterclockwise enclosing the $\pm  y_i^{1/2}, i=1, \ldots, N$
and $A$ is large enough so that the two contours $\Gamma$ and $\Upsilon=i\mbR+A$ do not cross each other.
\begin{figure}[htbp!]
\psfrag{a}{$-y_1^{1/2}$}
\psfrag{b}{$-y_2^{1/2}$}
\psfrag{h}{$-y_N^{1/2}$}
\psfrag{d}{$y_N^{1/2}$}
\psfrag{e}{$y_2^{1/2}$}
\psfrag{f}{$y_1^{1/2}$}
\psfrag{U}{$\Upsilon$}
 \includegraphics[height=6.5cm, width=11cm]{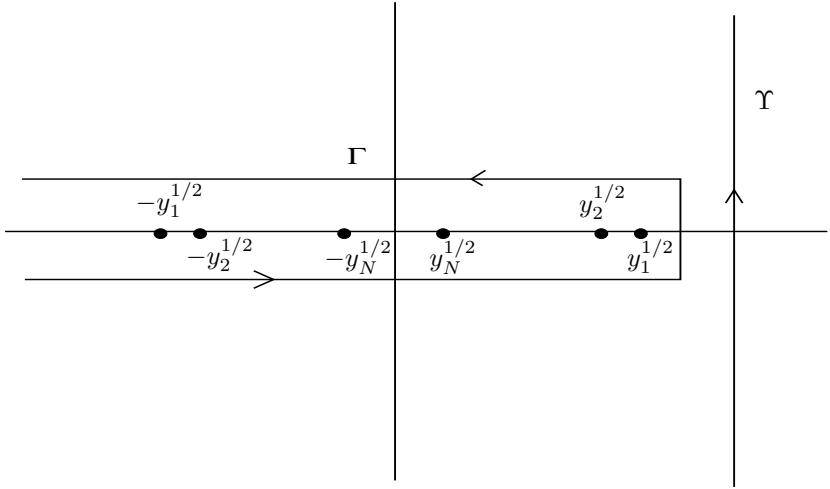}
\caption{The two contours defining the correlation kernel $K_N(u,v;H)$.}
\end{figure}

Due to the determinantal structure, correlation functions are not impacted by conjugation of the correlation kernel. In particular, for any $b \in \mathbb{R}$ and any integer $m\geq 1,$
$$\det \left ( K_N(u_i,u_j;H)\right)_{i,j=1}^m=\det \left (K_N(u_i,u_j;H)e^{\frac{2b(\sqrt{u_j}-\sqrt{u_i})}{S}}\right)_{i,j=1}^m.$$
In the sequel we consider the conjugation for some $b$ that will be defined in the asymptotic analysis.
We set 
\be
\label{def: conjkernel}
K_N^b(u_i,u_j;H):=K_N(u_i,u_j;H)e^{2b(\sqrt{u_j}-\sqrt{u_i})/S}.
\ee

The correlation kernel $K_N^b(u,v;H)$ only depends on $H$ through its spectral measure $\mu_N$. As in the case where $a$ is fixed (independent of $N$), the idea is to use the convergence of $\mu_N$ to the Mar\v{c}enko-Pastur distribution. In particular, one would like to make the replacement (outside a suitably negligible set)
$ \prod_{i=1}^N (w^2-y_i)=\exp\{N (1+o(1))\int \ln (w^2-y)d\rho_{\gamma}^{MP}(y)\}.$
In the next subsection, we prove that this replacement can be made in some sense provided $\Im (w^2)$ is not too small.  
\subsection{Concentration results \label{Sec: concentration}}

This subsection is devoted to the proof of the following Proposition \ref{Prop: concentration}, which precises the rate of convergence of $\mu_N$ to $\rho^{MP}$. Before exposing this Proposition we need a few definitions and notations.
Given a complex number $z=u+i \eta, $ $u \in \mathbb{R}, \eta >0$, we define the Stieltjes transform of $\mu_N$ by
\be \label{def: mN}
m_N(z):=\int \frac{1}{\lambda-z}d\mu_N(\lambda),\ee
and the Stieltjes transform of the limiting Mar\v{c}enko-Pastur distribution by 
\be \label{def: mMP}
m_{MP}(z):=\int \frac{1}{\lambda-z}d\rho^{MP}_{\gamma }(\lambda). \ee

Let $\epsilon >0$ be given (small enough).
\bp \label{Prop: concentration} Let $z=u+i\eta$ for some $u \in [u_-+\epsilon, u_+-\epsilon]$ and $\eta>0$. 
Then, there exist a constant $c_1$ and $C_o,C>0, c>0$ depending on $\epsilon $ only such that $\forall \delta<c_1\epsilon,$
$$\mathbb{P}\left (\sup_{u \in [u_-+\epsilon, u_+-\epsilon]}\Big | m_N(z)- m_{MP}(z) \Big | \geq \big(\delta+C_o\big|\frac{p}{N}-\gamma\big| \big)\right) \leq Ce^{-c\delta \sqrt{N \eta}},$$
for any $(\ln N)^4/N \leq \eta \leq 1.$
Furthermore, given $\eta \geq (\ln N)^4/N$, there exist constants $c>0,C>0$ and $K_o$ such that $\forall \kappa \geq K_o$, $$\mathbb{P}\left (\sup_{ |x|>(\epsilon/200)^2, y\geq \eta } |m_N(x+iy)|\geq \kappa \right) \leq Ce^{-c\sqrt{\kappa N \eta}}.$$\ep 

\paragraph{Proof of Proposition \ref{Prop: concentration}:}
To ease the reading, the proof is postponed to Appendix A. 
The proof is actually a modification of that of Theorems 3.1 and 4.1 in \cite{SchleinYau1}, where the Stieltjes transform of Hermitian Wigner matrices is considered. We thus simply indicate the main changes.

\subsection{Asymptotics of the correlation kernel when $\nu=p-N$ is bounded\label{sec: asanal}}
The aim of this subsection is here to prove the following Proposition.

Let $\epsilon>0$ be given (small) independent of $N$.
\bp \label{prop: asym_det} Let $u_* \in [u_-+\epsilon, u_+-\epsilon]=[\epsilon, 1-\epsilon]$ and $m\geq 1$ be given.
Consider a sequence $u=u_N$ such that $N^{1-\lambda}|u-u_*|=o(1).$
Then, there exist $b\in \mbR$, a set $\Omega_N$, and positive constants $C, c$ such that:\\ 
-the complement of $\Omega_N$ is negligible: $\mathbb{P}(\Omega_N^c)\leq Ce^{-cN^{\lambda/4}}$\\
-on $\Omega_N$ one has that
\begin{eqnarray*}&&\lim_{N \to \infty} \det \left( \frac{1}{N \rho^{MP}(u)}K_N^b\left (u+\frac{x_i}{N \rho^{MP}(u)},u+\frac{x_j}{N \rho^{MP}(u)};H\right)\right)_{i,j=1}^m \crcr
 &&=\det \left (\frac{\sin (\pi(x_i-x_j)}{\pi(x_i-x_j)} \right)_{i,j=1}^m.
\end{eqnarray*}
\ep 
\brem \label{rem: nana}In Remark \ref{rem: extu*} we explain how to extend Proposition \ref{prop: asym_det} to the case where $|u-u_*|=o(N^{-c})$ with $c>0$ arbitrarily small. This extension is needed for the proof of Theorem \ref{theo: spacing}.
\erem 
This subsection is devoted to the proof of Proposition \ref{prop: asym_det}. The proof is divided into two parts: first we obtain a new expression for the correlation kernel, which then allows to derive its asymptotics by a saddle point argument.
\subsubsection{Rewriting the kernel \label{sec: rewr}} 
Our strategy consists into two parts: we first replace the Bessel functions with their asymptotic expansion and then remove the singularity $1/(w-z)$ in the correlation kernel as it will be proved to prevent a direct saddle point analysis.
As this part is highly technical and needs a few notations, we summarize our main result in Lemma \ref{lem: resume} stated at the end of this subsection. The reader might skip this part and is simply referred to the above cited Lemma.

\begin{figure}[h]
\psfrag{a}{$-y_1^{1/2}$}
\psfrag{b}{$-y_2^{1/2}$}
\psfrag{h}{$-y_N^{1/2}$}
\psfrag{d}{$y_N^{1/2}$}
\psfrag{e}{$y_2^{1/2}$}
\psfrag{f}{$y_1^{1/2}$}
\psfrag{U}{$\Upsilon$}
\psfrag{V}{$\Gamma_c$}
\psfrag{m}{$x_1^+$}
\psfrag{n}{$x_1^-$}
\psfrag{o}{$x_0^-$}
\psfrag{p}{$x_0^+$}
 \includegraphics[height=6.5cm, width=11cm]{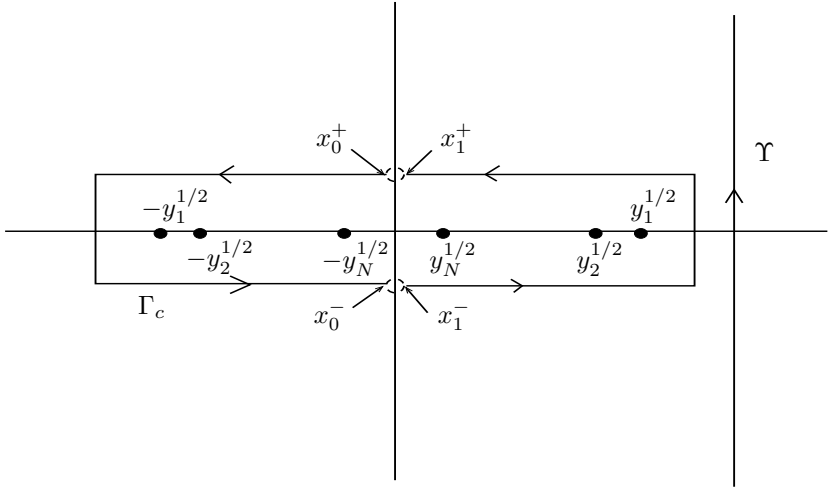}
\caption{The contour $\Gamma_c$ has been cut along the imaginary axis.}
\label{fig: cutGamma}
\end{figure}

Our first task to replace in (\ref{defKN}) Bessel functions with their asymptotic expansion given in Appendix \ref{app: B}, formula (\ref{asyBessel}).
This replacement can be made, up to a negligible error, provided the contour $\Gamma$ is cut in a small neighborhood of the imaginary axis: see Figure \ref{fig: cutGamma}. In the sequel (see Figure \ref{fig: cutGamma}) we call $x_o^{\pm}$ (resp. $x_1^{\pm}$) the endpoints of $\Gamma_c$ with $\Re (x_o)<0$ (resp. $\Re (x_1)>0$) with positive/negative imaginary part surrounding the imaginary axis: $x_o^{\pm}$ and $x_1^{\pm}$ will be chosen in subsection \ref{subsec: saddle}. 
We also call $\Gamma_1$ (resp. $\Gamma_2$) the part of $\Gamma_c$ lying to the right (resp. left) of $i\mbR$.
Then $\Gamma_1$ is the image of $\Gamma_2$ by $z\mapsto -z$ and $\Gamma_c=\Gamma_1\cup \Gamma_2.$
One obtains the following Lemma.

\bl \label{lem: Bessreplace}For any $b\in \mathbb{R},$ one has that 
\be\label{approxKN}
K_N^b(u,v;H)=(1+o(1)) \left ( K_N^{1,b}(u,v;H)+K_N^{2,b}(u,v;H)\right)
,\ee
where 
\begin{eqnarray}
&K_N^{1,b}(u,v;H)&=\frac{e^{\frac{2b(\sqrt v-\sqrt u)}{S}}}{4i^2\pi^{2} S}\int_{i\mbR+A}dw\int_{\Gamma_1}dz e^{\frac{-z^2+w^2}{S}}e^{-\frac{2w\sqrt v}{S}}e^{
\frac{2z\sqrt u}{S}}\crcr
&& \frac{1}{\sqrt{w}\sqrt{z} (uv)^{1/4}}\prod_{i=1}^N \frac{w^2-y_i}{z^2-y_i}\left (\frac{w}{z}\right)^{\nu} \frac{w+z}{w-z}, \crcr 
&K_N^{2,b}(u,v;H)&=\frac{e^{\frac{2b(\sqrt v-\sqrt u)}{S}}}{4i^2\pi^{2} S}\int_{i\mbR+A}dw\int_{\Gamma_1}dz e^{\frac{-z^2+w^2}{S}}e^{-\frac{2w\sqrt v}{S}}e^{-\frac{2z\sqrt u}{S}}\crcr
&& \frac{1}{\sqrt{w}\sqrt{z} (uv)^{1/4}}\prod_{i=1}^N \frac{w^2-y_i}{z^2-y_i}\left (\frac{w}{z}\right)^{\nu} \frac{w+z}{w-z}e^{\nu i \pi+ i \pi/2}. 
\end{eqnarray}
\el
 \paragraph{Proof of Lemma \ref{lem: Bessreplace}:}
Here we will call on some arguments already used in \cite{GBAPecheCPAM}, which we won't develop entirely, yet trying to be the most self-contained as possible.
Lemma \ref{lem: Bessreplace} essentially follows from the two claims exposed in the sequel and whose proof relies on the saddle point analysis in the next subsection.
\begin{claim} \label{claim1}
One has that
\begin{eqnarray}
&&\lim_{N \to \infty}\frac{e^{\frac{2b(\sqrt v-\sqrt u)}{S}}}{2i^2\pi^2 S^2}\int_{i\mbR+A}dw\int_{\Gamma\setminus \Gamma_c}dz e^{\frac{-z^2+w^2}{S}}K_{\nu}\Big (\frac{2w\sqrt v}{S}\Big )I_{\nu}\Big (\frac{2z\sqrt u}{S}\Big )\crcr
&& \prod_{i=1}^N \frac{w^2-y_i}{z^2-y_i}\left (\frac{w}{z}\right)^{\nu} \left (\frac{w+z}{w-z}-\frac{w-z}{w+z}\right)=0. 
\end{eqnarray}
\end{claim}
Claim \ref{claim1} essentially follows from the fact that $|I_{\nu}(z)|\leq e^{z}$ if $\Re (z)>0$ (and a similar bound for $K_{\nu}$)
and the fact that $\Gamma\setminus \Gamma_c$ lies far away from the critical points (a full justification can be derived 
from the saddle point analysis performed in Subsection \ref{subsec: saddle}, see Remark \ref{rem:asBes}). 

\begin{claim} \label{claim2}
One has that
\begin{eqnarray}
&&\frac{e^{\frac{2b(\sqrt v-\sqrt u)}{S}}}{2i^2\pi^2 S^2}\int_{i\mbR+A}dw\int_{ \Gamma_c}dz e^{\frac{-z^2+w^2}{S}}K_{\nu}\Big (\frac{2w\sqrt v}{S}\Big )I_{\nu}\Big (\frac{2z\sqrt u}{S}\Big )\crcr
&& \prod_{i=1}^N \frac{w^2-y_i}{z^2-y_i}\left (\frac{w}{z}\right)^{\nu} \left (\frac{w+z}{w-z}-\frac{w-z}{w+z}\right)\crcr
&&=\frac{e^{\frac{2b(\sqrt v-\sqrt u)}{S}}}{4i^2\pi^2 S}\int_{i\mbR+A}dw\int_{ \Gamma_1}dz e^{\frac{-z^2+w^2}{S}}e^{-\frac{2w\sqrt v}{S}+\frac{2z\sqrt u}{S}}\frac{1}{(uv)^{\frac{1}{4}}\sqrt w\sqrt z}\crcr
&& \prod_{i=1}^N \frac{w^2-y_i}{z^2-y_i}\left (\frac{w}{z}\right)^{\nu} \left (\frac{w+z}{w-z}-\frac{w-z}{w+z}\right)\left (1+O(N^{-\lambda/2})\right ).
\end{eqnarray}
\end{claim}
By a straightforward change of variables, we can reduce $\Gamma_c$ to $\Gamma_1$.  
Claim \ref{claim2} then follows by uniform asymptotic expansions of first the $w-$integral, using the asymptotics of the Bessel function $K_{\nu}$ given in (\ref{asyBessel}), and then of the $z-$integral. Again this follows from the saddle point analysis of Subsection \ref{subsec: saddle}. This saddle point analysis will prove that the correlation kernel does not vanish in the large $N$-limit. Then splitting $\Gamma_1$ into sufficiently many pieces lying to the left or to the right of the critical points and moving accordingly $\Upsilon=i\mbR +A$ to pass close enough to the same critical points yields that the error term is of order $N^{-\lambda/2}$. This yields Claim \ref{claim2} and Lemma \ref{lem: Bessreplace}.
We skip the detail. $\square$

\paragraph{} 
A saddle point analysis of the correlation kernel $K_N^{1,b}$ cannot be performed at that point, due to the singularity $1/(w-z)$. Indeed, assuming that $\prod_{i=1}^N (w^2-y_i)=e^{\int N\ln (w^2-y) d\rho^{MP}(y)(1+o(1))},$ it is not difficult to see that the $z-$ and $w-$integrands have the same critical points. Thus we first remove this singularity (see \cite{GBAPecheCPAM} for a more detailed explanation). \\
Define 
\begin{eqnarray}&&f^N_u(z):= z^2 -2\sqrt u z+ S\sum_{i=1}^N\ln (z^2-y_i), \, f_N(z)=f^N_u(z)+2b\sqrt{u}\crcr
&&f^N_v(z):= z^2 -2\sqrt v z+ S\sum_{i=1}^N\ln (z^2-y_i), \, \widetilde{f}_N(z)=f^N_v(z)+2b\sqrt{v}\crcr
&&g_{N,u}(w,z):=(w-b)\frac{(f_u^N)'(w)-(f_u^N)'(z)}{w-z}+(f_u^N)'(z),\crcr
 &&g_{N,u}^1(w,z)=\frac{\nu bS}{wz}+\frac{Sb}{2wz}\frac{w-z}{w+z},\:\theta(w,b)=-\frac{e^{\frac{2(w-b)(\sqrt u-\sqrt v)}{S}}-1}{2(w-b)(\sqrt u -\sqrt v)},\crcr
&&\theta^1(z,b)=\frac{e^{-\frac{2(z-b)(\sqrt u+\sqrt v)}{S}}-1}{2(z-b)(\sqrt u +\sqrt v)}.\label{def: ttesfonctions}\end{eqnarray}
\bl \label{lem: rewr} One has that
\begin{eqnarray}
&&K_N^{1,b}(u,v;H)=\crcr&&  \frac{1}{(uv)^{1/4}4i^2\pi^2S}\int_{i\mbR+A}dw\int_{\Gamma_1}dz e^{\frac{f_N(w)-f_N(z)}{S}}\crcr
&& \times \frac{w+z}{\sqrt{w}\sqrt{z }}\left (g_{N,u}(w,z)+g_{N,u}^1(w,z)\right )\left (\frac{w}{z}\right)^{\nu}\theta(w,b)\label{mainK11}\\
&&-\frac{(x_1^--b)}{4i^2\pi^2(uv)^{1/4}}\int_{i\mbR+A}dw\frac{w+x_1^-}{\sqrt{w}\sqrt{x_1^- }}\left (\frac{w}{x_1^-}\right)^{\nu}e^{\frac{f_N(w)-f_N(x_1^-)}{S}}\theta(w,b)\crcr
&&+\frac{(x_1^+-b)}{4i^2\pi^2(uv)^{1/4}}\int_{i\mbR+A}dw\frac{w+x_1^+}{\sqrt{w}\sqrt{x_1^+ }}\left (\frac{w}{x_1^+}\right)^{\nu}e^{\frac{f_N(w)-f_N(x_1^+)}{S}}\theta(w,b)
.\crcr
&&\label{KN1rewr}
\end{eqnarray}
Setting $\widetilde{K}_N^{2,b}(u,v;H):= e^{-\nu i \pi- i \pi/2}K_N^{2,b}(u,v;H),$ we have that
\begin{eqnarray}
&&\widetilde{K}_N^{2,b}(u,v;H)=\crcr && \frac{e^{-\frac{4b\sqrt u}{S}}}{4i^2\pi^2S}\int_{i\mbR+A}dw\int_{\Gamma_1}dz
\left (g_{N,v}(w,z)+g_{N,v}^1(w,z)\right )e^{\frac{\widetilde{f}_N(w)-\widetilde{f}_N(z)}{S}}\crcr
&& \times \frac{w+z}{\sqrt{w}\sqrt{z} (uv)^{1/4}}\left (\frac{w}{z}\right)^{\nu}\theta^1(z,b)\crcr
&&-\frac{e^{-\frac{4b\sqrt u}{S}}(x_1^--b)}{4i^2\pi^2(uv)^{1/4}}\int_{i\mbR+A}dw\frac{w+x_1^-}{\sqrt{w}\sqrt{x_1^- }}\left (\frac{w}{x_1^-}\right)^{\nu}e^{\frac{\widetilde{f}_N(w)-\widetilde{f}_N(x_1^-)}{S}}\theta^1(x_1^-,b)\crcr
&&+\frac{e^{-\frac{4b\sqrt u}{S}}(x_1^+-b)}{4i^2\pi^2(uv)^{1/4}}\int_{i\mbR+A}dw\frac{w+x_1^+}{\sqrt{w}\sqrt{x_1^+ }}e^{\frac{\widetilde{f}_N(w)-\widetilde{f}_N(x_1^+)}{S}}\left (\frac{w}{x_1^+}\right)^{\nu}\theta^1(x_1^+,b)
.\crcr
&&\label{KN2rewr}
\end{eqnarray}
\el 
\paragraph{Proof of Lemma \ref{lem: rewr}:}
We only consider $K_N^{1,b}$ (the arguments are similar for $K_N^{2,b}$).
To this aim we make the change of variables $w=b+\beta (w'-b)$, $z=b+\beta (z'-b)$ for some $\beta $ real close to $1$ and get the following : 
set 
\be E(w):= e^{\frac{w^2-2\sqrt vw}{S}}\prod_{i=1}^N (w^2-y_i);\,
G(z)=e^{\frac{z^2-2\sqrt uz}{S}}\prod_{i=1}^N (z^2-y_i).\ee
Then one has that
\begin{eqnarray}
&K_N^{1,b}(u,v;H)&= \frac{e^{\frac{2b(\sqrt v-\sqrt u)}{S}}}{(uv)^{\frac{1}{4}}4i^2\pi^{2}S}\int_{i\mbR+A}dw\int_{\Gamma_1}dz\frac{1}{\sqrt{wz }}\frac{E(w)}{G(z)}\left (\frac{w}{z}\right)^{\nu} \frac{w+z}{w-z}\crcr
&&=\frac{e^{\frac{2b(\sqrt v-\sqrt u)}{S}}}{4i^2\pi^{2}S(uv)^{\frac{1}{4}}} \int_{i\mbR+A}\int_{\Gamma_1}\frac{\beta dw dz}{\sqrt{(b+\beta(w'-b))(b+\beta(z'-b)) }}\crcr
&&\times\left (\frac{b+\beta(w'-b)}{b+\beta(z'-b)}\right)^{\nu} \frac{2b+\beta(w'+z'-2b)}{w'-z'}
\frac{E(b+\beta(w'-b))}{G(b+\beta(z'-b))}\crcr
&&-\frac{e^{\frac{2b(\sqrt v-\sqrt u)}{S}}}{4i^2\pi^{2}S(uv)^{\frac{1}{4}}}\int_{i\mbR+A}\int_{x_1^+}^{b+\beta (x_1^+-b)} \frac{dw dz}{\sqrt{wz }}\left (\frac{w}{z}\right)^{\nu} \frac{w+z}{w-z}
\frac{E(w)}{G(z)}\crcr
&&+\frac{e^{\frac{2b(\sqrt v-\sqrt u)}{S}}}{4i^2\pi^{2}S(uv)^{\frac{1}{4}}}\int_{i\mbR+A}\int_{x_1^-}^{b+\beta (x_1^--b)} \frac{dw dz}{\sqrt{wz }}\left (\frac{w}{z}\right)^{\nu} \frac{w+z}{w-z}
\frac{E(w)}{G(z)}.\crcr
&&
\end{eqnarray}
Differentiating with respect to $\beta$ (close to $1$), we find that (see \cite{Johansson} Section 2 for the detail)
\begin{eqnarray}
&& (uv)^{1/4}K_N^{1,b}+(\sqrt v-\sqrt u)\frac{d}{d\sqrt v}\Big ( (uv)^{1/4}K_N^{1,b}\Big )=
\crcr
&&\frac{e^{\frac{2b(\sqrt v-\sqrt u)}{S}}}{4i^2\pi^{2}S^2}\int_{i\mbR+A}dw\int_{\Gamma_1}dz\frac{1}{\sqrt{wz}}\left (\frac{w}{z}\right)^{\nu} \frac{w+z}{w-z}
\frac{E(w)}{G( z)}\crcr
&&\times \Bigl (\Big [-2(w+z-b)+2\sqrt u+S\sum_{i=1}^N \frac{2y_i(w+z)-b(2wz+y_i)}{(w^2-y_i)(z^2-y_i)}\Big ]\label{derivsec}\\
&& -\frac{\nu bS}{wz}-\frac{Sb}{2wz}\frac{w-z}{w+z}\Bigr) 
\crcr 
&&-\frac{(x_1^+-b)e^{\frac{2b(\sqrt v-\sqrt u)}{S}}}{4i^2\pi^{2}S}\int_{i\mbR+A}dw \frac{1}{\sqrt{wx_1^+ }}\left (\frac{w}{x_1^+}\right)^{\nu} \frac{w+x_1^+}{w-x_1^+}
\frac{E(w)}{G(x_1^+)}\crcr
&&+\frac{(x_1^--b)e^{\frac{2b(\sqrt v-\sqrt u)}{S}}}{4i^2\pi^{2}S}\int_{i\mbR+A}dw\frac{1}{\sqrt{wx_1^- }}\left (\frac{w}{x_1^-}\right)^{\nu} \frac{w+x_1^-}{w-x_1^-}
\frac{E(w)}{G(x_1^-)}.
\end{eqnarray}

By the definition of $f_u^N$, the term in the bracket (\ref{derivsec}) can also be written 
\begin{eqnarray}&(\ref{derivsec})=&- \frac{(w-b)(f_u^N)'(w)-(z-b)(f_u^N)'(z) }{w-z}\crcr
&&=-(w-b)\frac{(f_u^N)'(w)-(f_u^N)'(z)}{w-z}-(f_u^N)'(z). 
\end{eqnarray}

Integrating by parts yields Lemma \ref{lem: rewr}.$\square$

To summarize this highly technical subsection, we have proved that 
\bl \label{lem: resume} Assume that $\nu$ is a bounded integer. Then
$$K_N^b(u,v;H)=(1+o(1)) \left ( K_N^{1,b}(u,v;H)+K_N^{2,b}(u,v;H)\right),$$
where the kernels $K_N^{1,b}$ and $K_N^{2,b}$ are defined in Lemma \ref{lem: rewr}.
\el 

\subsubsection{Saddle point analysis of the correlation kernel. \label{subsec: saddle}}
We are now in position to perform the saddle point analysis of the correlation kernel $K_N^b(u,v;H)$
and prove Proposition \ref{prop: asym_det}. Let $\tau$ be a given real number independent of $N$ and 
assume that $v-u=\tau/(N \rho^{MP} (u))$. We mainly focus on the asymptotics of the first term (\ref{mainK11}) in $K_N^{1,b}$. The asymptotic expansion of the two other terms (\ref{KN1rewr}) in $K_N^{1,b}$ and of $K_N^{2,b}$ will be an easy corrollary of the arguments used hereafter.

\paragraph{}We first examine the saddle point analysis for the approximate exponential term :
\be \label{def: f}f(w)=f_u(w):= w^2-2(w-b)\sqrt u +a^2 \int \ln (w^2-y) d\rho^{MP}(y).\ee
Note that $f$ is the almost sure limit of the exponential term arising in the definition of (\ref{mainK11}) and should thus (approximately) lead its asymptotic analysis.
Then using the fact that, when $\gamma=1$, a Mar\v{c}enko-Pastur random variable has the same law as a squared semi-circle random variable, one gets that
\be \label{rewr: f}f'(w):=2w-2\sqrt u+ a^2\int \frac{2w}{w^2-y}d\rho^{MP}(y)=2w-2\sqrt u +4a^2(w- \sqrt{w^2-1}).\ee
Given a point $u \in [\epsilon, 1-\epsilon], \epsilon >0$
it is easy to see that $f$ admits two conjugate critical points:
$$w_c^{\pm}=\frac{(1+2a^2)\sqrt u \pm 2a^2i \sqrt{1+4a^2-u}}{1+4a^2}=\sqrt u \pm a^2i \pi \sqrt u \rho^{MP}(u) +O(a^2).$$
More precisely, one has that $\Im (w_c^{\pm})= \pm a^2 \pi \sqrt u \rho^{MP}(u)+O(a^4),$
where the $O$ is uniform in the bulk (depending on $\epsilon $ only).
One can also check that $$f''(w_c^{\pm})= 2(1+2a^2)+\frac{4a^2w_c^{\pm}}{\sqrt{(w_c^{\pm})^2-1}}=2+O(a^2),$$
where the $O(\cdot)$ is uniform as long as $u\in [\epsilon, 1-\epsilon], \epsilon >0$ and depends on $\epsilon $ only.\\
We now give the relevant contours for a saddle point analysis of the approximate exponential terms.
Let $\Upsilon=\Re (w_c^+)+it, t\in \mathbb{R}$ and $\Gamma_1=\Gamma_1^+\cup \overline{\Gamma_1^+}$ with 
$\Gamma_1^+=\{w_{c}^{+}( t), t \in [\frac{\epsilon}{2}, 1-\frac{\epsilon}{2}]\}\cup \{w_{c}^{+}({1-\frac{\epsilon}{2}})+x, x>0\}\cup \{i \Im (w_{c}^{+}({\frac{\epsilon}{2}}))+x, \epsilon/180 \leq x\leq \Re (w_{c}^{+}({\frac{\epsilon}{2}}))\}$. 
Then it is easy to check that $\Re f$ achieves its maximum (resp. minimum) at $w_c^{\pm}$ on $\Upsilon$ (resp. $\Gamma$).
 
\paragraph{}We now turn to the saddle point analysis for the true exponential term: 
\be \label{def: fN}f_N(w):=w^2-2(w-b)\sqrt u +\frac{a^2}{N}\sum_{i=1}^N \ln (w^2-y_i).\ee
One has that $f_N'(w)=2(w-\sqrt u)-2wa^2m_N(w^2).$
Using the concentration results of Proposition \ref{Prop: concentration}, we now show
that the first derivatives of $f_N$ and $f$ are close (on a suitable set).  
Let $\zeta >0$ be a very small number and set 
\be S_{\epsilon, \zeta}=\Bigl \{z=x+iy,\, x\in [u_-+\frac{\epsilon}{180}, u_+-\frac{\epsilon}{8}],\, \zeta a^2\leq y \leq 1\Bigr \}, \label{def: edk}\ee
with $u_-=0$ and $u_+=1$ here as $\gamma=1$.
Define for $H=H_N=\frac{WW^*}{N}$,
\begin{eqnarray}&\Omega_N=\Big \{\!\!\!\!\! &|m_N(z)-m(z)|\leq N^{\frac{-\lambda}{4}},\,\text{if}\, \frac{\epsilon}{180}\leq \Re (z)\leq 1-\frac{\epsilon}{8},\, \zeta a^2\epsilon \leq \Im (z)\leq 1;\crcr
&& \!\!\!\sup y_i(H) \leq K;\, \sup_{ |x|>(\frac{\epsilon}{200})^2, y\geq \zeta a^2 } |m_N(x+iy)| \leq K \Big \}.\label{def: OmegaN}
\end{eqnarray}
Using Proposition \ref{Prop: concentration} and Lemma 7.3 in \cite{SchleinYau0}, we can deduce that there exist $K$ large enough and constants $c,C>0$ such that 
$$\mathbb{P}(\Omega_N^c)\leq Ce^{-cN^{\lambda/4}}.$$
From now on we assume that $H_N \in \Omega_N.$ As a consequence, there exists $C>0$ such that for any $w \in S_{\epsilon, \zeta}$, 
\be |f_N'(w)-f'(w)|\leq \frac{Ca^2}{N^{\lambda/4}}.\label{contrini}\ee
By Cauchy's formula, one can then deduce that for any integer $l\geq 1$
\be |f_N^{(l)}(w)-f^{(l)}(w)|\leq \frac{C}{a^{2(l-2)}N^{\lambda/4}}\label{contrderi} .\ee
Here the value of the constant $C$ varies from line to line.

\paragraph{}Let us now show that $f_N$ admits two conjugate critical points $w_{c,N}^{\pm}$ which are very close to $w_c^{\pm}.$
The proof directly follows the arguments of \cite{SchleinYau3}, Section 3.
The critical points of $f_N$ are the solutions of the fixed point equation:
$$w=F_N(w)=\sqrt u -\frac{a^2}{N}\sum_{i=1}^N \frac{w}{w^2-y_i}.$$
This equation clearly admits $2N-1$ real solutions which are interlaced with $-\sqrt y_N, \ldots, -\sqrt y_1, \sqrt y_1, \ldots, \sqrt y_N.$ It admits also two non real solutions $w_{c,N}^{\pm}$. We now show that $w_{c,N}^{\pm}$ are very close to $w_c^{\pm}.$
Set $F(w)=\sqrt u -a^2\int \frac{w}{w^2-y}d\rho^{MP}(y).$
Define $\Theta=\{ z\in \mathbb{C}, \, |\Re (z) -\sqrt u|\leq C_0 a^2,\: \zeta a^2\leq \Im (z)\leq C_0a^2\},$ for some large constant $C_0$.
Then, on $\Theta,$ $|F_N-F|\leq Ca^2N^{-\lambda/4}$.
As on $\Theta$, $\Re \left (F(z)\right)=\sqrt u +O(a^2)$ and $\Im (F(z))=a^2\pi\sqrt u\rho^{MP}(u)+o(a^2)$ we deduce that $F_N(\Theta)\subset \Theta.$
Now it is an easy consequence of (\ref{contrderi}) that $F_N$ is a contraction: one has that $|F_N'(w)-F'(w)|\leq CN^{-\lambda/4}$.
Thus $F_N$ admits a unique fixed point with $\Im (w_{c,N}^{+})>\zeta a^2$. Furthermore, it is an easy fact that
$$|w_{c,N}^{+}-w_c^+|\leq \frac{Ca^2}{N^{\lambda/4}}.$$
 
\paragraph{}We now slightly modify the contours $\Gamma$ and $\Upsilon$ for the saddle point analysis of $f_N$.
Set $\Upsilon_N=\Big\{\Re (w_{c,N}^{+})+it, t\in \mathbb{R}\Big \}$. Then $\Upsilon_N$ lies within a $C^1-$distance of at most $Ca^2N^{-\lambda/4}$ from $\Upsilon$. 
Furthermore, setting $w=\Re (w_{c,N}^{+})+it$,
$$\Re \left(\frac{d}{dt}f_N(w)\right) =-t\Big \{1-\frac{a^2}{N}\sum_{i=1}^N \left (\frac{1}{|w-\sqrt y_i|^2}+\frac{1}{|w+\sqrt y_i|^2}\right)\Big \}.$$
By the monotonicity of $t\mapsto \sum_{i=1}^N \left (\frac{1}{|w-\sqrt y_i|^2}+\frac{1}{|w+\sqrt y_i|^2}\right),$ we conclude that
$\Re (f_N)$ is maximum on $\Upsilon_N$ at $w_{c,N}^{\pm}$.
More precisely, using that $\Re f_N''(w)\geq 1$ for any $w\in \Upsilon \cap \{\Im (w)\geq \zeta a^2\}$, one has that, 
if $t$ stands for $\Im(w)$, 
\begin{eqnarray}&\Re (f_N(w)-f_N(w_{c,N}^+))&=-\Re \left( \int_{0}^{t-\Im (w_{c,N}^+)}sds \int_{0}^1 f_N''(w_{c,N}^++ius) du\right)\crcr
&&\leq -\frac{(t-\Im (w_{c,N}^+))^2}{4}. \label{maincontri2}
\end{eqnarray}
This implies in particular that 
$$\Re \left( f_N\left (\Re (w_{c,N}^{+})+i\zeta a^2\right )-f_N(w_{c,N}^{+})\right )\leq -\left (\zeta -a^{-2}\Im (w_{c,N}^{+})\right )^2a^4/4.$$
Let then $w\in \Upsilon \cap \{z, 0\leq \Im (z) \leq \zeta a^2\}$.Then, with a slight abuse of notation when $\Im(w)=0$,
$$\Big |\exp{\Big\{f_N(w) -f_N\left (\Re (w_{c,N}^{+})+i\zeta a^2\right )\Big\}}\Big |\leq e^{\zeta^2 a^4}.$$
If one chooses $\zeta$ in (\ref{def: edk}) small enough so that for any $v \in [\epsilon, 1-\epsilon], $ $\zeta a^2\leq \Im (w_{c,N}^{+})/8$, one can then deduce that the contribution of the contour $\Upsilon \cap \{z, \Im (z) \leq \zeta a^2\}$ is negligible.
Thus, using (\ref{maincontri2}), the main contribution to the $w-$integral comes from a neighborhood of width $\sqrt{S}$ of the critical point.

We now turn to the $z-$contour. Let 
\begin{eqnarray}&\Gamma_N^+=&\Big \{w_{c,N}^{+}(t), \epsilon/2\leq t \leq 1-\epsilon/2\Big \}
\cup \Big \{w_{c,N}^{+}\left ({1-\epsilon/2}\right )+x, x>0\Big \}\crcr&&\cup \Big \{i \Im \Big (w_{c,N}^{+}\left ({\epsilon/2}\right )\Big )+x, \epsilon/180 \leq x\leq \Re \Big (w_{c,N}^{+}\left ({\epsilon/2}\right )\Big )\Big \}.\end{eqnarray}  
We can now define $x_1^{\pm}:=\epsilon/180\pm i\Im \left (w_{c,N}^+\Big({\epsilon/2}\Big)\right ).$ The choice of 180 here is arbitrary: we only need a sufficiently large constant. We also set $x_o^{\pm}=-x_1^{\mp}.$

By construction, for any $t, \, \frac{\epsilon}{2}\leq t \leq 1-\frac{\epsilon}{2},$
$$\Re \left( -f_N(w_{c,N}^{+}(t))+f_N\Big (w_{c,N}^{+}\Big )\right)\leq -c(\sqrt u -\sqrt{t})^2,$$ for some constant $c>0$ small enough. This follows from the fact that for any $t, \, \frac{\epsilon}{2}\leq t \leq 1-\frac{\epsilon}{2},$ $\Im (w_{c,N}^{+}(t))\geq \zeta a^2$ and $\Re \left( \frac{dw_{c,N}^{+}(t)}{dt}\right)=\frac{1+O(a^2)}{2\sqrt{t}}.$ 
Moreover $\Re \left(f_N(w_{c,N}^{+}(1-\frac{\epsilon}{2}))-f_N\big (w_{c,N}^+(1-\frac{\epsilon}{2})+x\big )\right)\leq -cx^2,$ for some constant $c>0$ small enough. This follows from the fact that on $\Omega_N$ $a^2 zm_N(z^2)\to 0$ uniformly along $\Gamma_N$ as $m_N$ is bounded.
Thus the contribution of $\{w_{c,N}^{+}({1-\frac{\epsilon}{2}})+x, x>0\}$ is negligible as $N \to \infty.$ 
For the same reason, the contribution of $\{i \Im (w_{c,N}^{+}({\frac{\epsilon}{2}}))+x, \epsilon/180 \leq x\leq \Re (w_{c,N}^+({\frac{\epsilon}{2}}))\}$ is negligible in the large $N$ limit. 
As for the $w-$integral, the main contribution to the $z-$integral comes from a neighborhood of width $\sqrt{S}$ of the critical point.
\brem \label{rem:asBes} This analysis justifies Claim \ref{claim1} as $\Re \left( f_N(w_{c,N}^+)-f_N(x_1^+)\right)\leq -c$ for some $c>0.$
One can choose $\Gamma\setminus \Gamma_c$
as $\{x_1^++ia^2t, 0<t<A-\Im (x_1^+)\}\cup\{x_o^++ia^2t, 0<t<A-\Im (x_1^+)\}\cup \{x+ia^2A, \Re (x_o^+)\leq x\leq \Re (x_1^+)\}$ plus its conjugate, for some large enough constant $A$. There exists $C>0$ such that $\Re \left( f_N(x_1^+)-f_N(z)\right) \leq Ca^2, \,\forall z\in \Gamma \setminus \Gamma_c$, thus yielding a negligible contribution.
\erem

\paragraph{}We can now conclude to the asymptotic expansion of the correlation kernel and prove Proposition \ref{prop: asym_det}. Let $u_*$ be as in Proposition \ref{prop: asym_det}.\\
We first consider the asymptotics for the first term (\ref{mainK11}) in $K_N^{1,b}$. 
We now fix $b$ as follows:
$$b=\Re \left ( w_{c,N}^{+}(u_*)\right).$$
Thanks to this choice, $|\Re ( w_{c, N}^+)-b|=o(N^{\lambda-1}).$
Therefore the function $\theta$ has no impact on the saddle point argument exposed in the above, neither do the functions $g_{N,u}$ and $g_{N,u}^1$.  At the critical points, one deduces from (\ref{def: ttesfonctions}) that:
\begin{eqnarray*}
&&g_{N,u}(w_{c,N}^{\pm}, w_{c,N}^{\pm})=(w_{c,N}^{\pm}-b)f_N''(w_{c,N}^{\pm})\simeq i\Im \left(w_{c,N}^{\pm} \right)f_N''(w_{c,N}^{\pm});\crcr
&& g_{N,u}(w_{c,N}^{\pm}, w_{c,N}^{ \mp})=0;\crcr
&&g_{N,u}^1(w_{c,N}^{\pm}, w_{c,N}^{\pm})=g_{N,u}^1(w_{c,N}^{\pm}, w_{c,N}^{ \mp})=O(S)=O(\Im \left(w_{c,N}^{\pm} \right)/N);\crcr
&&\theta (w_{c,N}^{\pm},b)=\frac{e^{i \frac{2(\sqrt u-\sqrt v) \Im (w_{c,N}^{\pm})}{S}}-1}{2 i\Im \left ( w_c^{N , \pm}\right )(\sqrt u -\sqrt v)}(1+o(1)). 
\end{eqnarray*}
Let us now consider the 4 combined contributions of the different critical points. 
The contribution to (\ref{mainK11}) from $z=w=w_{c,N}^{\pm}$ gives (at the leading order) 
\be \label{combcontr}\frac{\pm i}{4i^2\pi^2S}\left ( \frac{\sqrt{2\pi S}}{\sqrt{f''(w_{c,N}^{\pm})}}\right)^2 
\frac{1}{(uv)^{1/4}}f''(w_{c,N}^{\pm}) \frac{e^{\frac{2 \Im w_{c,N}^{\pm}(\sqrt u-\sqrt v)}{S}}-1}{(\sqrt u -\sqrt v)}
.\ee
The contribution to (\ref{mainK11}) from $z=\overline{w}=w_{c,N}^{\pm}$ is $O(N^{1-\lambda/2})$ due to the fact that $g_{N,u}$ anneals at that point.
Combining the above yields that:
$$\lim_{N \to \infty}
\dfrac{1}{N \rho^{MP}(u)}(\ref{mainK11})=\lim_{N \to \infty}\frac{1}{N \rho^{MP}(u)}(\ref{combcontr})=\frac{\sin (\pi \tau)}{\pi \tau }.  
$$

The contributions of the two other terms in (\ref{KN1rewr}) is exponentially small since there exists a constant $c>0$ depending on $\epsilon$ such that:
$$\Re \left (f_N(w_{c,N}^{\pm}(\frac{\epsilon}{2})\right )-\Re \left ( f_N (w_{c,N}^{+})\right) >c.$$
This finishes the asymptotic expansion of $K_N^{1,b}.$\\
Let us now turn to the asymptotics for $K_N^{2,b}.$
The function $\theta^1$ is not bounded. Nevertheless there exists $\chi>0$ depending on $\epsilon$ such that 
$\Big | e^{-\frac{4b\sqrt u}{S}}\theta^1(z,b)\Big | \leq e^{-\chi S^{-1}}.$
This follows from the fact that $\Re (b-z)<b$ along the contour $\Gamma_N$. One can now 
use the same saddle point arguments as for the study of $K_N^{1,b}$ to show that 
$$\frac{1}{N \rho^{MP}(u)}K_N^{2,b}(u,v;H)= O(e^{-\chi S^{-1}}).$$
This is enough to ensure Proposition \ref{prop: asym_det}. $\square$

\brem \label{rem: extu*}In the case where $|u-u_*|=o(N^{-c})$ for some $c\leq 1-\lambda$, we choose a sequence 
$\hat{u}(N)$ such that $N^{1-\lambda}|u-\hat{u}(N)|\to 0,$ $|u_*-\hat{u}(N)|=o(N^{-c})$ and $\hat{u}(N)\in [u_-+\epsilon, u_+-\epsilon]$ for all $N$. We can then replace $u_* \to \hat{u}(N)$ in all the asymptotic analysis (in particular $b=b(N)=\Re \left(w_{c, N}^+(\hat{u}(N))\right)$). This has no impact on the validity of the arguments.
\erem
\subsection{The case where $\nu\to \infty$\label{sec: gennu}}
We now turn to the case where $\nu$ is unbounded (that is either $\gamma=1$ and $\nu \to \infty$ or $\gamma\not=1$) and study the asymptotics of the correlation kernel (\ref{defKN}). 
The aim of this section is to prove the following Proposition. Let $\epsilon>0$ be given (small).
\bp \label{prop: sinegamma}Assume that $\gamma\geq 1$ and $\lambda<1/2$. Let $u_*\in [u_-+\epsilon, u_+-\epsilon]$ and $u=u(N)$ be a sequence such that $\lim_{N \to \infty}N^{1-\lambda}(u-u_*)=0.$ 
Then there exists $b \in \mbR$ such that outside a set of negligible probability,
$$\lim_{N\to \infty}\dfrac{1}{N \rho_{\gamma}^{MP}(u)}K_N^b(u,u+\frac{\tau}{N \rho_{\gamma}^{MP}(u)};H)
=\frac{\sin (\pi \tau)}{\pi \tau}.$$
\ep 
\brem Proposition \ref{prop: sinegamma} extends to the case where $|u-u_*|=o(N^{-c})$ for some $c\leq 1-\lambda$ by the same arguments as in Remark \ref{rem: extu*}. 
\erem

The whole subsection is devoted to the proof of Proposition \ref{prop: sinegamma}. The proof follows the same steps as in subsection \ref{subsec: saddle} and we explain the main changes only.
\subsubsection{Rewriting the kernel}
Again the basic argument is to replace Bessel functions with their asymptotic expansion (large order large argument) and then derive the asymptotics of the correlation kernel by a saddle point argument.\\
We again cut the contour $\Gamma$ in a small neighborhood of the imaginary axis as on Figure \ref{fig: cutGamma}. We also call
$\Gamma_c$ this cut contour and the endpoints $x_1^{\pm}, x_o^{\pm}$ will be defined in the sequel. Let us denote again $\Gamma_1$ the part of $\Gamma_c$ lying to the right of the imaginary axis, so that $\Gamma_c:=\Gamma_1\cup (-\Gamma_1)$.
We now consider the uniform asymptotic expansion (\ref{asyBessel2}) of Appendix B.
We assume for a while that we can replace the modified Bessel functions with their asymptotic expansion in the correlation kernel (this will be proved in Lemma \ref{Lem: asyB2} below).
Thus, setting $v=u+\frac{\tau}{N \rho_{\gamma}^{MP}(u)}$, we consider the kernel:
\begin{eqnarray}
& \widetilde{K_N}(u,v;H):=& \frac{1}{4i^2\pi^2 S(uv)^{1/4}}\int_{i \mbR+A}dw \int_{\Gamma_1}dz \frac{e^{\frac{w^2-z^2+2(z-b)\sqrt{u}-2(w-b)\sqrt{v}}{S}}}{\sqrt w\sqrt z}\crcr
&& e^{\frac{(\gamma-1)^2 N^{\lambda}}{4u_*}(1/w-1/z)} \prod_{i=1}^N \frac{w^2-y_i}{z^2-y_i}\left(\frac{w}{z} \right)^{\nu}
\frac{w+z}{w-z}\label{KN1bis}\\
&&-\frac{1}{4i^2\pi^2 S(uv)^{1/4}}\int_{i \mbR+A}dw \int_{\Gamma_1}dz
\frac{e^{\frac{w^2-z^2+2(z-b)\sqrt{u}-2(w-b)\sqrt{v}}{S}}}{\sqrt w\sqrt z}\crcr
&& e^{\frac{(\gamma-1)^2 N^{\lambda}}{4u_*}(1/w-1/z)} \prod_{i=1}^N \frac{w^2-y_i}{z^2-y_i}\left(\frac{w}{z} \right)^{\nu}
\frac{w-z}{w+z}.\label{KN2bis}\end{eqnarray}
We call again $K_N^{1,b}$ (resp. $K_N^{2,b}$) the kernel given by (\ref{KN1bis}) (resp. \ref{KN2bis}). Set 
$$\Psi(w, z)=\exp{\{\frac{(\gamma-1)^2 N^{\lambda}}{4u_*}(1/w-1/z)\}}.$$
We can now slightly modify the arguments of subsection \ref{sec: rewr} to deduce that:
\begin{eqnarray}
&&K_N^{1,b}(u,v;H)=\crcr&&  \frac{1}{(uv)^{1/4}4i^2\pi^2S}\int_{i\mbR+A}dw\int_{\Gamma_1}dz e^{\frac{f_N(w)-f_N(z)}{S}}\Psi(w,z)\crcr
&& \times \frac{w+z}{\sqrt{w}\sqrt{z }}\left (g_{N,u}(w,z)+g_{N,u}^2(w,z)\right )\left (\frac{w}{z}\right)^{\nu}\theta(w,b)\label{mainK11bis}\\
&&-\frac{(x_1^--b)}{4i^2\pi^2(uv)^{1/4}}\int_{i\mbR+A}dw\frac{w+x_1^-}{\sqrt{w}\sqrt{x_1^- }}\left (\frac{w}{x_1^-}\right)^{\nu}e^{\frac{f_N(w)-f_N(x_1^-)}{S}}\theta(w,b)\Psi(w,x_1^-)\crcr
&&+\frac{(x_1^+-b)}{4i^2\pi^2(uv)^{1/4}}\int_{i\mbR+A}dw\frac{w+x_1^+}{\sqrt{w}\sqrt{x_1^+ }}\left (\frac{w}{x_1^+}\right)^{\nu}e^{\frac{f_N(w)-f_N(x_1^+)}{S}}\theta(w,b)\Psi(w,x_1^+)
,\crcr
&&\label{KN1rewrbis}
\end{eqnarray}
where $$g_{N,u}^2(w,z)=g_{N,u}^1(w,z)+\frac{SN^{\lambda}(\gamma-1)^2}{4u_*}\left (\frac{1}{wz}-b\frac{w+z}{z^2w^2}\right).$$
We leave $K^{2,b}_N$ unchanged since the singularity $1/(w+z)$ will not prevent its direct saddle point analysis.

\subsubsection{Saddle point analysis}
We shall now perform the saddle point analysis of the correlation kernels hereabove. The arguments follow closely those of Subsection \ref{sec: asanal}.
The approximate exponential term to be considered here is given by 
$$h_u(z):= z^2-2\sqrt u (z-b) +a^2 \int \ln (z^2-y) d\rho_{\gamma}^{MP}(y)+a^2(\gamma-1)\ln z.$$
We are indeed going to show that the perturbative term coming from $\Psi$ does not play a role in the asymptotics.
Let $m_{\gamma}(z)$ be the Stieltjes transform of the Mar\v{c}enko-Pastur distribution $\rho_{\gamma}^{MP}$.
One can check that
$$m_{\gamma}(z)=-2+\frac{\gamma-1}{2z}+2\frac{\sqrt{(z-(1+\sqrt \gamma)^2/4)(z-(1-\sqrt \gamma)^2/4)}}{z}.$$
The function $h_u$ admits two non real critical points which are conjugate :
\be
w_c^{\pm}=w_c^{\pm}( u)=\sqrt u(1-2a^2)\pm ia^2 \pi \sqrt u \rho^{MP}_{\gamma}(u)+O(a^4).\label{def: critgamma}\ee  
We now define the contours relevant for the saddle point analysis of the approximate exponential term.
Let $\Upsilon=\Re (w_c^+)+it, t \in \mbR$ be oriented positively from bottom to top.
The contour $\Gamma_1$ is defined as in Section \ref{subsec: saddle}: we set $\Gamma_1=\Gamma_{1}^+\cup \overline{\Gamma_1^+}$ where $\Gamma_1$ is oriented counterclockwise and
\begin{eqnarray*}&\Gamma_1^+=&\Big\{ w_c^{+}( t), t\in [u_-+\epsilon/90, u_+-\epsilon/2]\Big \}\cup \Big \{w_c^+(u_+-\epsilon/2)+x, x>0\Big \} \crcr&&\cup \Big \{i \Im \Big (w_{c}^{+}\left (u_-+{\frac{\epsilon}{90}}\right )\Big )+x, \epsilon/180 \leq x\leq \Re \Big (w_{c}^{+}\left ({u_-+\frac{\epsilon}{90}}\right )\Big )\Big \}.\end{eqnarray*}
We first consider the $z$-contour. There exists $c(\epsilon) >0$ depending on $\epsilon$ only such that for all points $t\in [u_-+\epsilon/90, u_+-\epsilon/2]$ in the bulk of the Mar\v{c}enko-Pastur distribution $\Re\left ( h_u''(w_c^+( t))\right)>c(\epsilon)$.
Using the fact that  
\be
\frac{d}{dt}\Re \left( h_u(w_c^+( t)) \right)=(\sqrt{t} -\sqrt u) \Re \frac{1}{\sqrt{t} h_u''(w_c^+( t))}, 
\ee 
$-\Re (h_u)$ achieves its maximum at $z=w_c^{\pm}$ along the first part of $\Gamma_1^+.$
There now remains to show that $-\Re (h_u)$ also decreases along $\Big \{w_c^+(u_+-\epsilon/2)+x, x>0\Big \}$ and 
$\Big \{w_c^+(u_-+\frac{\epsilon}{90})-x, 0<x<\Re \Big (w_{c}^{+}\left ({u_-+\frac{\epsilon}{90}}\right )\Big )-\frac{\epsilon}{180}\Big \}.$
This assertion follows from the fact that $m_{\gamma}$ is bounded along this curve.  
Using the same arguments as in Section \ref{subsec: saddle}, it is easily seen that $\Re (h_u)$ admits its maximum on $\Upsilon$ at $w_c^{\pm}.$ 
All the above reasoning holds unchanged if $\gamma$ is replaced by $p/N$. Thus, without loss of generality, we assume in the rest of this section that $p/N=\gamma.$

\paragraph{}We now consider the true exponential term. We set
$$h_{N,u}(z)=h_{N}(z):=z^2-2\sqrt u (z-b) + \frac{a^2}{N}\sum_{i=1}^N\ln (z^2-y_i)+a^2(\gamma-1)\ln z.$$
Again we can show that, provided $H$ belongs to a suitable set, the derivatives of $h_u$ and $h_{N,u}$ are very close and that
$h_{N,u}$ admits two non real critical points $w_{c,N}^{\pm}$ very close to $w_{c}^{\pm}$.
Let $\zeta>0$ be given (small) and $S_{\epsilon, \zeta}$ be given by (\ref{def: edk}).
Define for $H=H_N=\frac{WW^*}{N}$,
\begin{eqnarray}&\Omega_{N, \gamma}=\Big \{& \sup y_i(H) \leq K;\, \sup_{ |x|>(\frac{\epsilon}{200})^2, y\geq \zeta a^2} |m_N(x+iy)| \leq K;\crcr
&& |m_N(z)-m(z)|\leq N^{\frac{-\lambda}{4}},\, \forall\, z\text{ with }\zeta a^2\epsilon\leq \Im (z)\leq 1\crcr
&& \text{ and }u_-+\frac{\epsilon}{180} \leq \Re (z)\leq u_+-\frac{\epsilon}{8} \Big \}.\label{def: OmegaNbis}
\end{eqnarray}
Using Proposition \ref{Prop: concentration} and Lemma 7.3 in \cite{SchleinYau00}, we can deduce that there exist $K$ large enough and constants $c,C>0$ such that 
$$\mathbb{P}(\Omega_{N, \gamma}^c)\leq Ce^{-cN^{\lambda/4}}.$$
From now on we assume that $H_N \in \Omega_{N, \gamma}.$
Mutatis mutandis, all the arguments used to consider the true exponential term when $\nu$ is bounded (starting with (\ref{contrini})) can be copied to consider the true exponential term when $\nu$ is not bounded.
In particular, one has that 
$$|w_{c,N}^{\pm}-w_{c}^{\pm}|\leq \frac{Ca^2}{N^{\lambda/4}}, \,|h_{N,u}^{(l)}(w)-h_u^{(l)}(w)|\leq \frac{C}{a^{2(l-2)}N^{\lambda/4}}
\, \forall l\geq 1, \forall w \in S_{\epsilon,\zeta}.$$
We now consider the subsequent modified contours.
We set $\Gamma_N=\Gamma_N^+\cup \overline{\Gamma_N^+}$ with
\begin{eqnarray}&\Gamma_N^+=&\Big \{w_{c,N}^{+}(t), u_-+\frac{\epsilon}{90}\leq t \leq u_+-\frac{\epsilon}{2}\Big \}
\cup \Big \{w_{c,N}^{+}\Big ({u_+-\frac{\epsilon}{2}}\Big )+x, x>0\Big \}\crcr&&\cup \Big \{i \Im \left  (w_{c,N}^{+}\Big (u_-+{\frac{\epsilon}{90}}\Big )\right )+x, \epsilon/180 \leq x\leq \Re \left (w_{c,N}^{+}\Big (u_-+{\frac{\epsilon}{90}}\Big )\right )\Big \};\crcr
&\Upsilon_N=&\Re \left (w_{c,N}^{+}\right )+it, t\in \mbR.\end{eqnarray}
Accordingly, we now choose $x_1^{\pm}=\pm i \Im \Big (w_{c,N}^{+}\left (u_-+{\frac{\epsilon}{90}}\right )\Big )+\frac{\epsilon}{180}$ and $x_o^{\pm}=-x_1^{\mp}.$ \\  
Then the main contribution to the exponential term of (\ref{mainK11bis}) comes from a neighborhood of width $\sqrt{S}$ of $w_{c, N}^{\pm}$ for both the $w-$ and $z-$integral.
Thus one is left with showing that the function $\Psi$ does not impact on the saddle point analysis. Let $\mathcal{O}_w=\{w'\in \mbC, |w'-w|\leq N^{\alpha}\sqrt{S}\}$ for some $\alpha$ that we determine hereafter.
In order to ensure that $\Psi$ has no impact on the saddle point analysis outside $\mathcal{O}_{w_{c,N}^{\pm}}$, it is enough that $\alpha >\lambda/2$.
Now one can check that 
$$\forall \, w,z \in \mathcal{O}_{w_{c,N}^{\pm}}, \Big |\Psi(w,z)-1\Big |\leq CN^{3\lambda/2+\alpha-1}<<1,$$
if $\alpha$ and $\lambda$ are small enough (in particular $\lambda<1/2$). This finishes the proof that 
the main contribution to (\ref{mainK11bis}) comes from a neighborhood of width $\sqrt{S}$ of $w_{c, N}^{\pm}$ for both the $w-$ and $z-$integral.

We now set $b=\Re \left(w_{c, N}^+(u_*)\right),$ so that $\Big |b-w_{c, N}^+\Big|=o(N^{\lambda-1}),$ which implies in particular that the function $\theta$ does not impact the saddle point argument.
Let us now consider the contribution of the functions $g_{N,u}(w,z)+g_{N,u}^2(w,z)$ and $\theta(w,z)$ close to the critical points.
One has that 
\begin{eqnarray*}
&g_{N,u}(w_{c,N}^{\pm}, w_{c,N}^{\pm})+g_N^2(w_{c,N}^{\pm}, w_{c,N}^{\pm})&= \left(w_{c,N}^{\pm} -b\right)h_N''(w_{c,N}^{\pm})+O(S)\crcr
&&+\frac{SN^{\lambda}(\gamma-1)^2}{4u_*} \frac{w_{c,N}^{\pm}-2b}{(w_{c,N}^{\pm})^3}\crcr
&&=\left(w_{c,N}^{\pm} -b\right)h_N''(w_{c,N}^{\pm})+O(SN^{\lambda});\crcr
&g_{N,u}(w_{c,N}^{\pm}, w_{c,N}^{\mp})+g_N^2(w_{c,N}^{\pm}, w_{c,N}^{\mp})&=O(SN^{\lambda});\crcr
&\theta (w_{c,N}^{\pm},b)&=\frac{e^{i \frac{2(\sqrt u-\sqrt v) \Im (w_{c,N}^{\pm})}{S}}-1}{2 \left(w_{c,N}^{\pm} -b\right)(\sqrt u -\sqrt v)}(1+o(1)).\end{eqnarray*}

We can copy the end of the proof of Proposition \ref{prop: asym_det} to show that on $\Omega_{N,\gamma}$ 
$$\lim_{N \to \infty}\dfrac{1}{N \rho_{\gamma}^{MP}(u)} K^{1,b}_N(u,v;H)
=\dfrac{\sin (\pi \tau)}{\pi \tau}.$$

To study the asymptotics of $K^{2,b}_N$, we use the same saddle point analysis argument. Note that this is possible since
$\Upsilon_N\cap (-\Gamma_N)=\emptyset$.
We consider as leading exponential term for both the $w$- and $z$-integrals the function $h_{N,u}$ as above. Note that the $w$-integrand has then a ``perturbative term'' $e^{2(w-b)(\sqrt{u}-\sqrt{v})/S}$ which does not play a role in the saddle point analysis.
It is easy to deduce that 
$$\lim_{N \to \infty}\dfrac{1}{N \rho_{\gamma}^{MP}(u)} K^{2,b}_N(u,v;H)=0.$$
This follows from the fact that the contribution from equal critical points is negligible (due to the $(w-z)$ factor), while that of conjugate critical points is of order $1/N$ due to the rescaling and the fact that 
$\Big |e^{S^{-1}(h_{N,u}(z_{c, N}^{\pm}(u))-h_{N,v}(w_{c,N}^{\mp}(v))}\Big |\leq C$ for some $C>0$ independent of $N$.

\begin{figure}[htbp!]
\psfrag{U}{$\Upsilon$}
 \includegraphics[height=6cm, width=11cm]{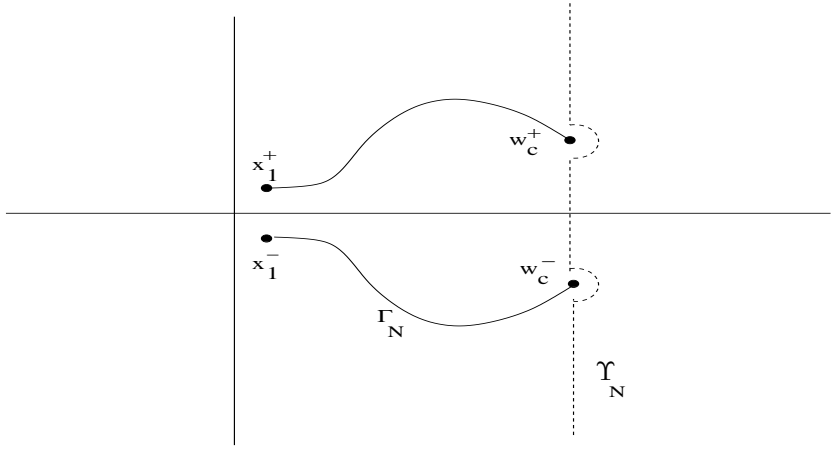}
\caption{The contour $\Upsilon_N$ has been deformed around $w_{c,N}^{\pm}.$}
\label{fig:deform}
\end{figure}

\paragraph{}Thus we are now left with showing the following Lemma. 
\bl \label{Lem: asyB2} One has that
$$\lim_{N \to \infty}\dfrac{1}{N \rho_{\gamma}^{MP}(u)}\Bigg| \widetilde{K_N}(u,v;H) -K_N^b(u,v;H)\Bigg|=0.$$
\el
\paragraph{Proof of Lemma \ref{Lem: asyB2}:}Note that it is enough to show that $\big| \widetilde{K_N}(u,v;H) -K_N^b(u,v;H)\big|=o(N)$.
First Claim \ref{claim1} can be translated with no modification to the case where $\nu $ is unbounded. Also we can reduce $\Gamma_c$ to $\Gamma_1$ in $K_N^b$ using the change of variables $z\mapsto -z$.
There remains to prove the counterpart of Claim \ref{claim2}.
Set \begin{eqnarray*}&Z=\dfrac{2\sqrt u z}{(\gamma-1)a^2}, \,Z_*=\dfrac{2\sqrt u_* z}{(\gamma-1)a^2} &W=\dfrac{2\sqrt v w}{(\gamma-1)a^2}, \,W_*=\frac{2\sqrt u_* w}{(\gamma-1)a^2}\crcr
&A_1'(Z):=\dfrac{e^{\nu (Z-\frac{1}{2Z})}}{\sqrt{2\pi \nu Z}}, &A_1(Z)=\frac{e^{\nu (Z-\frac{1}{2Z_*})}}{\sqrt{2\pi \nu Z}}
\crcr
&A_2'(W):=\sqrt{\pi}\dfrac{e^{-\nu (W-\frac{1}{2W})}}{\sqrt{2 \nu W}}, &A_2(W)=\sqrt{\pi}\frac{e^{-\nu (W-\frac{1}{2W_*})}}{\sqrt{2 \nu W}}.\end{eqnarray*}
In $\widetilde{K_N}$ we have replaced the Bessel functions $I_{\nu}$ and $K_{\nu}$ with the ``approximations'' 
$A_1$ and $A_2$. It is an easy computation that, along $\Upsilon_N$ and $\Gamma_N$, one has:
\be \label{mojerr}\Big |\frac{I_{\nu}(\nu Z)}{A'_1(Z)}-1\Big |\leq O(S), 
\Big |\frac{K_{\nu}(\nu W)}{A'_2(W)}-1\Big |\leq O(S).\ee 
We first consider the contribution of $\Gamma_{N}^l:=\Gamma_N\cap\{ \Re(x_1)\leq \Re(z) \leq \Re(w_{c,N}^+)\}$ to $\widetilde{K_N}(u,v;H) -K_N^b(u,v;H)$ where the new contour $\Upsilon_N^{\epsilon'}$ is slightly deformed around $w_{c,N}^{\pm}$ so that $d(\Upsilon_N, \Gamma_N^l) =\epsilon'\sqrt{S}$ for some small $\epsilon'$ (see Figure \ref{fig:deform}).
It is then easy to deduce from (\ref{mojerr}) and the previous saddle point analysis that 
\begin{eqnarray}
&\Bigg |\dfrac{e^{\frac{2b(\sqrt v -\sqrt u)}{S}}}{2i^2\pi^2 S^2}\displaystyle{\int_{\Upsilon_N^{\epsilon'}}dw\int_{\Gamma_N^l}dz} &e^{\frac{-z^2+w^2}{S}}\prod_{i=1}^N \frac{w^2-y_i}{z^2-y_i}\left (\frac{w}{z}\right)^{\nu} \frac{w+z}{w-z}\crcr
&& \left ( A'_2(W)I_{\nu}\Big (\frac{2z\sqrt u}{S}\Big )-A'_2(W)A_1'(Z)\right) \Bigg|\crcr
&&\leq O(\sqrt S).\label{(1)}\end{eqnarray} 
Indeed replacing $A_1$ and/or $A_2$ with $A'_1$ and/or $A'_2$ does not impact on the saddle point argument.
Similarly
\begin{eqnarray}
&\Bigg |\dfrac{e^{\frac{2b(\sqrt v -\sqrt u)}{S}}}{2i^2\pi^2 S^2}\displaystyle{\int_{\Upsilon_N^{\epsilon'}}dw\int_{\Gamma_N^l}dz} &e^{\frac{-z^2+w^2}{S}}\prod_{i=1}^N \frac{w^2-y_i}{z^2-y_i}\left (\frac{w}{z}\right)^{\nu} \frac{w+z}{w-z}\crcr
&& \!\!\left( A'_2(W)I_{\nu}\Big (
\frac{2z\sqrt u}{S}\Big )- K_{\nu}\Big (\frac{2w\sqrt v}{S}\Big )I_{\nu}\Big (
\frac{2z\sqrt u}{S}\Big )\right) \Bigg |\crcr
&&\leq  O(\sqrt{S}),\label{(2)}\end{eqnarray}
by using (\ref{mojerr}), the fact that $I_{\nu}(\nu Z)/A'_1(Z) $ is uniformly bounded for $N$ large enough and the previous saddle point analysis. 
Lastly 
\begin{eqnarray}
&\Bigg |
\dfrac{e^{\frac{2b(\sqrt v -\sqrt u)}{S}}}{2i^2\pi^2 S^2}\displaystyle{\int_{\Upsilon_N^{\epsilon'}}dw\int_{\Gamma_N^l}dz} &e^{\frac{-z^2+w^2}{S}}\left (A_2(W)A_1(Z)-A'_2(W)A_1'(Z)\right )\crcr
&& \prod_{i=1}^N \frac{w^2-y_i}{z^2-y_i}\left (\frac{w}{z}\right)^{\nu}  \frac{w+z}{w-z}\Bigg|\crcr
&&\!\!\!\!\leq O\left (\frac{N^{\lambda+\alpha} \sqrt{S}|u-u_*|}{\epsilon' \sqrt S}+\frac{N^{\lambda}|u-u_*|}{N^{\alpha }\sqrt S}+ N^{\lambda/2}\right)\crcr &&\label{majerenfin}\\
&&\!\!\!\!\leq O\left (N^{\alpha+2\lambda-1}+N^{\frac{3\lambda}{2}-\alpha}+N^{\frac{\lambda}{2}}\right)=o(N).\label{(3)}
\end{eqnarray}
In (\ref{majerenfin}) we have separated the cases where $d(w,z)\geq N^{\alpha}\sqrt S$ or not, where $z\in \Gamma_N^l$ and $w\in \Upsilon_N^{\epsilon'}$ and used the fact that $|u-v|\simeq N^{-1}$ and $|u-u_*|=o(N^{\lambda-1}).$
Replacing $\frac{w+z}{w-z}$ with $\frac{w-z}{w+z}$ in (\ref{(1)}), (\ref{(2)}) and (\ref{(3)}) yields similar estimates.

We now turn to the contribution of $\Gamma_{N}^r:=\Gamma_N\cap\{ \Re(z) \geq \Re(w_{c,N}^+\}$ to $\widetilde{K_N}(u,v;H) -K_N^b(u,v;H)$. Note that $\Upsilon_N$ can be moved to the left of $\Re(w_{c,N}^+)$ up to adding the residue at $w=z.$
Define $\Upsilon_N^{-\epsilon'}$ as the contour obtained by reflecting $\Upsilon_N^{\epsilon'}$ with respect to the line $\Re(z)=\Re (w_{c,N}^+).$
Then 
\begin{eqnarray}
&\displaystyle{\int_{A+i\mbR}dw\int_{\Gamma_N^r}dz} &e^{\frac{-z^2+w^2}{S}}\prod_{i=1}^N \frac{w^2-y_i}{z^2-y_i} \left (\frac{w}{z}\right)^{\nu} \frac{w+z}{w-z}\crcr
&&\left (K_{\nu}\Big (\frac{2w\sqrt v}{S}\Big )I_{\nu}\Big (
\frac{2z\sqrt u}{S}\Big )-A'_1(Z)A'_2(W)\right )\crcr 
&=\displaystyle{\int_{\Upsilon_N^{-\epsilon'}}dw\int_{\Gamma_N^r}dz} &e^{\frac{-z^2+w^2}{S}}\prod_{i=1}^N \frac{w^2-y_i}{z^2-y_i}\left (\frac{w}{z}\right)^{\nu} \frac{w+z}{w-z}\crcr
&&\left (K_{\nu}\Big (\frac{2w\sqrt v}{S}\Big )I_{\nu}\Big (
\frac{2z\sqrt u}{S}\Big )-A'_1(Z)A'_2(W)\right ) \crcr
&&\label{termepareil}\\
&+2i\pi \displaystyle{\int_{\Gamma_N^r}dz}2z \Big (K_{\nu}\Big (\frac{2z\sqrt v}{S}\Big )&I_{\nu}\Big (
\frac{2z\sqrt u}{S}\Big )-A'_1(Z)A'_2\Big (\frac{2z\sqrt v}{a^2 (\gamma-1)}\Big )\Big )\label{termeapart}.
\end{eqnarray}
The analysis of (\ref{termepareil}) is similar to that of $\Gamma_N^l\cap \Upsilon_N^{\epsilon'}.$
We thus have to consider (\ref{termeapart}).
One has that \be \Big |(\ref{termeapart})\Big |=O\left ( S\int_{\Gamma_N^r}dz2z A'_1(Z)A'_2\Big (\frac{2z\sqrt v}{a^2 (\gamma-1)}\Big )\right )\label{termeapart2}.\ee
As the integrand in the r.h.s of (\ref{termeapart2}) has no singularity, one can thus move $\Gamma_N^r$ to the line joining the two critical points $w_{c,N}^{\pm}.$
Using now that $|b- \Re (w_{c, N}^+)|=o(N^{\lambda-1})$,
one can easily see that there exists a constant $C>0$ such that
$$e^{\frac{2b(\sqrt v -\sqrt u)}{S}}(\ref{termeapart2})S^{-2}\leq C.$$
One gets similar estimates replacing $\frac{w+z}{w-z}$ with $\frac{w-z}{w+z}$.
This yields Lemma \ref{Lem: asyB2}. $\square$
\section{Proof of Theorem \ref{theo: uni2} and Theorem \ref{theo: spacing} \label{Sec: pfth}}

We only give the proof of Theorem \ref{theo: uni2} when $\gamma=1$ and $\nu$ is bounded. The extension to other parameters $\gamma$ is straightforward.
The proof of universality for the spacing distribution can easily be deduced from \cite{Johansson} and \cite{SchleinYau3} and is based on the extension of Proposition \ref{prop: asym_det} given in Remark \ref{rem: nana}.
Let then $f \in L^{\infty}(\mathbb{R}^2)$ with compact support and $S_N^{(2)}(f, u, \widetilde{\rho_N})$ be
defined by $(\ref{les})$ with $u\in [\epsilon, 1-\epsilon]$ and $\widetilde{\rho_N}=N \rho^{MP}(u).$
Then if $\mbE$ denotes the expectation w.r.t. the distribution of $\widetilde{M}_N$ in (\ref{defensem}),
$$\mathbb{E}S_N^{(2)}(f, u , \widetilde{\rho_N}) 
=\int dP_N(H)\int_{\mathbb{R}^N_+}f_N^H(x_1, \ldots, x_N)S_N^{(2)}(f, u, \widetilde{\rho_N})[x] dx^N,$$
where $f_N^H(\cdot)$ (resp. $S_N^{(2)}(f, u, \widetilde{\rho_N})[\cdot]$) is the density function (\ref{densdefwish}) (resp. local eigenvalue statistic) of the deformed Wishart ensemble.
Then, using (\ref{def: OmegaN}) and setting $$d\overline{P}_N(H)=\frac{1}{Z_N} dP_N(H)1_{H \in \Omega_N},$$ where $Z_N$ is the normalizing constant, we get that $\exists\, C>0$ such that
\begin{eqnarray}&&\Big |\mathbb{E}S_N^{(2)}(f,u, \widetilde{\rho_N})-
\int d \overline{P}_N(H)\int_{\mbR_+^N}f_N^H(x)S_N^{(2)}(f,u, \widetilde{\rho_N})[x]dx^N \Big | \crcr
&& \leq CN^2 |f| _{\infty}e^{-cN^{\lambda/4}}=o(1).\label{alpha}
\end{eqnarray}
Now, \begin{eqnarray}&&\lim_{N \to \infty}\int d \overline{P}_N(H)\int_{\mathbb{R}^N_+}f_N^H(x)S_N^{(2)}(f,u, \widetilde{\rho_N})[x]dx^N\crcr
&&=\lim_{N \to \infty}\int d\overline{P}_N(H)\int_{\mathbb{R}^2}\frac{f(t_1,t_2)}{\widetilde{\rho_N}^2}R_N^{(2)}\left(u+\frac{t_1}{\widetilde{\rho_N}},u+\frac{t_2}{\widetilde{\rho_N}}; H\right )dt_1 dt_2 \crcr
&&=\int_{\mathbb{R}^2}f(t_1,t_2)\left (1-\left (\frac{\sin \pi(t_1-t_2)}{ \pi(t_1-t_2)}\right)^2\right) dt_1dt_2
  \label{beta}   \end{eqnarray}
where we used the definition of correlation functions,
Proposition \ref{prop: asym_det} and the fact that $f$ has compact support.  
\paragraph{}There now remains to extend the result to non Gauss divisible ensembles in order to prove the full Theorem \ref{theo: uni2}. The argument exactly follows the arguments of \cite{SchleinYau3}. We recall them for seek of completeness.
(\ref{alpha}) and (\ref{beta}) show that the sine kernel holds for the complex measure $e^{t\mathcal{L}}G_td\mu^{\otimes Np}\otimes e^{t\mathcal{L}}G_td\mu^{\otimes Np}$ if $t=N^{-1+\lambda}$. More precisely, let $p_{N,t}(x)$
denote the density function of the eigenvalues $x=(x_1, \ldots, x_N)$
w.r.t. this measure and let $R_{N,t}^{(2)}$ be its two point correlation function. Similarly, we define $p_{N,c}(x)$ (resp. $p_N(x)$) and  $R_{N,c}^{(2)}$ (resp. $R_{N,F}^{(2)}$) for the eigenvalue density and two point correlation function w.r.t. truncated complex measure $F_{c, N,p}=f_{c, N,p}d\mu^{\otimes Np}\otimes f_{c, N,p}d\mu^{\otimes Np}$ (resp. w.r.t. measure $F_{N,p}=F^{\otimes Np}$).
We also set in the following
$$u(t_1)=u+\frac{t_1}{\rho_N} \text{ and } u(t_2)=u+\frac{t_2}{\rho_N } \text{ with }\rho_N=N\rho^{MP}_{1,1}(u).$$
Then
$$\Bigg|\int \Big[ R^{(2)}_{N,F}\Big(u(t_1),u(t_2)\Big)-R^{(2)}_{N,t}\Big(u(t_1),u(t_2)\Big)\Big]f(t_1,t_2)dt_1 dt_2\Bigg| 
\le (I)+ (II),$$
where
$$(I):=\Bigg|\int \Big[ R^{(2)}_{N,F}\Big(u(t_1), u(t_2) \Big)-R^{(2)}_{N,c}\Big(u(t_1), u(t_2)\Big)\Big]f(t_1,t_2)dt_1dt_2 \Bigg|,$$
$$(II):=\int \Big| R^{(2)}_{N,c}\Big(u(t_1),u(t_2)\Big) -R^{(2)}_{N,t}\Big(u(t_1), u(t_2)\Big)\Big|\, |f(t_1,t_2)|dt_1dt_2.$$
It is easy to see from Proposition \ref{prop: approxF} that
$$(I)\le N^2 \|f\|_\infty 2Ce^{-cN^c} \le C'e^{-c'N^c}$$
with some $C', c'>0$ as $N\to\infty$.
To estimate $(II)$, we use the fact that 
\begin{eqnarray}
&(II)^2 \le &   \int  \Big[ \frac{R_{N,c}^{(2)}}{R_{N,t}^{(2)}}\Big( u(t_1), u(t_2)\Big)-1\Big]^2  
R_{N,t}^{(2)}\Big( u(t_1), u(t_2)\Big) |f(t_1, t_2)| d t_1 d t_2 \Big] \crcr
&&\label{long1}\\
&& \times \Big[ \int R_{N,t}^{(2)}\Big( u(t_1), u(t_2)\Big) |f(t_1,t_2)| d t_1 d t_2 \Big].\label{long}
\end{eqnarray}
The factor (\ref{long}) is bounded using (\ref{alpha}) and (\ref{beta}).
We now use the beautiful idea of \cite{SchleinYau3} to bound (\ref{long1}). Set $\varrho:=\rho^{MP}_{1,1}(u).$
As $f$ is bounded, one has that
\begin{eqnarray}\label{NN}
 &(\ref{long1})&\leq  C   \Bigg[N^2\varrho^2 \int  
  \Big[ \frac{R^{(2)}_{N,c}(z,y)}{R^{(2)}_{N,t}(z,y)}-1 \Big ]^2R^{(2)}_{N,t}(z,y)dz dy \Bigg]^{1/2}\crcr
&&\leq C \Bigg [N^2\varrho^2 \int \left (\frac{p_{N,c}(x)}{p_{N,t}(x)} -1\right) ^2 p_{N,t}(x)dx\Bigg ]^{1/2}\crcr
&&\leq C\Bigg [N^2\varrho^2 \int \frac{|e^{t\mathcal{L}}G_t -f_{c,N,p}|^2}{e^{t\mathcal{L}}G_t}d\mu^{\otimes Np}\Bigg ]^{1/2}\crcr
&&\leq CN^{-1+4\lambda},
\end{eqnarray}
where $C>0$ is a constant that varies from line to line. Here the basic argument is that the distance $D(f,g):=\int |f/g-1|^2 g$ between two probability measures $f$ and $g$ decreases when taking marginals as well as when passing from the matrix ensemble to the induced joint eigenvalue density. Finally, we used the estimate of Proposition \ref{prop: approxF}. This completes the proof of Theorem \ref{theo: uni2}. $\square$

\appendix
\section{Appendix : Proof of Proposition \ref{Prop: concentration}}
To ease the reading, we here recall Proposition \ref{Prop: concentration}.
\bp Let $z=u+i\eta$ for some $u \in [u_-+\epsilon, u_+-\epsilon]$ and $\eta>0$. 
Then, there exist a constant $c_1$ and $C_o,C>0, c>0$ depending on $\epsilon $ only such that $\forall \delta<c_1\epsilon,$
$$\mathbb{P}\left (\sup_{u \in [u_-+\epsilon, u_+-\epsilon]}\Big | m_N(z)- m_{MP}(z) \Big | \geq \big(\delta+C_o\big|\frac{p}{N}-\gamma\big| \big)\right) \leq Ce^{-c\delta \sqrt{N \eta}},$$
for any $(\ln N)^4/N \leq \eta \leq 1.$
Furthermore, given $\eta \geq (\ln N)^4/N$, there exist constants $c>0,C>0$ and $K_o$ such that $\forall \kappa \geq K_o$, $$\mathbb{P}\left (\sup_{ |x|>(\epsilon/200)^2, y\geq \eta } |m_N(x+iy)|\geq \kappa \right) \leq Ce^{-c\sqrt{\kappa N \eta}}.$$
\ep 

\paragraph{Proof: }The proof follows closely that of Theorems 3.1, 4.1 and 4.6 in \cite{SchleinYau1}. We only prove here the first statement. The second one is simple adaptation of the latter statement and ideas given in the proof of Theorem 4.6 in \cite{SchleinYau1} (see also the proof of Theorem 2.1 in \cite{SchleinYau0}).\\ 
Using a discretization scheme with step $\delta/4\eta^2$ (see \cite{SchleinYau1}), it is possible to show that it is actually enough to prove that 
$$\mathbb{P}\left (\Big | m_N(z)- m_{MP}(z) \Big | \geq \big (\delta+C_o \big|\frac{p}{N}-\gamma\big| \big)\right) \leq Ce^{-c\delta \sqrt{N \eta}},$$
for any $(\ln N)^4/N \leq \eta \leq 1$ and for a given $z=u+i\eta$ where $u\in [u_-+\epsilon, u_+-\epsilon].$
Denote by $C_k$ the $k$th column of $W/\sqrt N$, $C_k:=\frac{1}{\sqrt N}(W)_k,$ then one has that:
\be\label{ide}H_N=\sum_{k=1}^p C_kC_k^*.\ee 
We define also $\mathbf{R}_N(z)=(H_N-zI)^{-1}$ and for any integer $k=1, \ldots, p,$ $\mathbf{R}_N^{(k)}(z):=(H_N -C_kC_k^*-zI)^{-1}.$
Then $m_N(z)= \frac{1}{N}\tTr \mathbf{R}_N(z)$ and one has that
\be \label{eq: resol}1+zm_{N}(z)=\frac{p}{N}-\frac{1}{N}\sum_{k=1}^p \frac{1}{1+C_k^*\mathbf{R}_N^{(k)}(z)C_k}.\ee
Equation (\ref{eq: resol}) simply follows from (\ref{ide}) and the identity $\mathbf{R}_N(z)(H_N-zI)=I.$
In addition we denote by $y_1^{(k)}\geq y_2^{(k)}\geq \cdots \geq y_N^{(k)}$ the ordered eigenvalues of $H_N-C_kC_k^*$ and set $\mu_N^{(k)}:=\frac{1}{N}\sum_{i=1}^N\delta_{y_i^{(k)}}.$
By the well-known interlacing property of eigenvalues, for any $k=1, \ldots, p,$
$y_1\geq y_1^{(k)}\geq y_2 \geq y_2^{(k)}\geq \cdots \geq y_N \geq y_N^{(k)}.$
Denote by $F_N$ (resp. $F_N^{(k)}$) the p.d.f. of the spectral measure $\mu_N$ (resp. $\mu_N^{(k)}$). Then one has that
\be \label{err}\Big |N F_N(x)-NF_N^{(k)}(x)\Big |\leq 1, \quad \forall x \in \mathbb{R}.\ee
Last call $v_{i}^{(k)}, i=1, \ldots, N$ a set of orthonormal eigenvectors associated to the ordered eigenvalues
$y_i^{(k)}$ and define 
$$\xi_i^{(k)}:= |<v_i^{(k)}, \sqrt{N} C_k>|^2.$$
From (\ref{eq: resol}), one gets that  
\be \label{eq: resol2}
1+zm_{N}(z)=\frac{p}{N}-\frac{1}{N}\sum_{k=1}^p \frac{1}{1+\frac{1}{N}\sum_{i=1}^N\frac{\xi_i^{(k)}}{y_i^{(k)}-z}}.
\ee
The latter formula (\ref{eq: resol2}) is the counterpart of formula (2.6) in \cite{SchleinYau1}.
The proof of \cite{SchleinYau1} can be summarized into 2 basic ideas.\\
First the random vectors $C_k$ are centered with i.i.d. entries and they are independent of $\mathbf{R}_N^{(k)}$. It is known that the random variables
$C_k^*\mathbf{R}_N^{(k)}(z)C_k$ concentrate around their means which is given by $ \mbE |W_{11}|^2\tTr \mathbf{R}_N^{(k)}/N$. The speed of concentration is explicited under assumption $({\bf A_1})$ in Lemma 4.2 in \cite{SchleinYau1} for Wigner matrices. The extension to Wishart matrices states as follows. 
\bl \label{lem:conc}Set $X^{(k)}=X=\frac{1}{N}\sum_{k=1}^p \frac{\sigma^{-2}\xi_i^{(k)}-1}{y_i^{(k)}-z}$ where $z=u+i\eta, u\in [u_-+\epsilon, u_+-\epsilon].$ Then there exists a positive constant $c$ (depending on $\epsilon$) so that for every
$\delta>0$ we have
$$\mathbb{P}[|X| \geq \delta] \leq 5 e^{-c\min\{\delta \sqrt{N \eta/\gamma}, \delta^2N \eta/\gamma\}} $$
if $N\eta\geq (\ln N)^2$ and $N$ is sufficiently large (independently of $\delta$).
\el 
The proof of Lemma \ref{lem:conc} is postponed to the end of Appendix A.
\brem \label{rem:iid}Lemma \ref{lem:conc} is established under the assumption that the real and imaginary parts of the components of $C_k$ are i.i.d. (in addition to the gaussian decay assumption $({\bf A_1})$). This is the reason for our assumption $(H_3).$
\erem
One can then use (\ref{err}) to show that for any $k=1, \ldots, p,$ $|\tTr \mathbf{R}_N^{(k)}-\tTr \mathbf{R}_N|\leq \eta^{-1}.$ 
By Lemma \ref{lem:conc}, $|X^{(k)}|\leq \delta$, $\forall k=1, \ldots,p$ in (\ref{eq: resol2}) with high probability. A bootstrap argument exposed in \cite{SchleinYau00}, Section 2, yields that with high probability, the Stieltjes transform $m_N(z)$ satisfies 
\be \label{stiel+err}1+zm_{N}(z)=\frac{p}{N}-\frac{p}{N}\frac{1}{1+\sigma^2m_N(z)}+\Delta,\ee
where $|\Delta|\leq C'\delta $ for some $C'>0$ is a small error term and $\sigma^2=\mbE |W_{11}|^2 =1/4$.\\
The second basic idea is the stability of the equation (\ref{stiel+err}). 
The equation
\be \label{stiellim}1+zm(z)=\frac{p}{N}-\frac{p}{N}\frac{1}{1+\sigma^2m(z)}\ee admits a unique solution  
satisfying $\Im (m(z))>0$ whenever $\Im (z)>0$. This solution is $m_{MP, p/N}(z)$, that is the Stieltjes transform of the Mar\v{c}enko-Pastur distribution with parameter $p/N$. Now, the stability of equation (\ref{stiel+err}) implies that there exists a constant $C$ such that for any $z\in \{ z=u+i\eta, u\in [u_-+\epsilon, u_+-\epsilon], (\ln N)^4/ N \leq \eta \leq 1\},$
$$ |m_N(z)-m_{MP, p/N}(z)|\leq C \Delta.$$ The constant $C$ here depends on $\epsilon $ only.
Now, as $|z-u_{\pm}| \geq \epsilon$ and $1/|z|=O(1)$, there exists $C_o>0$ depending on $\epsilon $ only such that $$|m_{MP, p/N}(z)-m^{MP}(z)|\leq C \big|p/N-\gamma\big|.$$ 
This finishes the proof of Proposition \ref{Prop: concentration} provided we show Lemma \ref{lem:conc}. $\square$

\paragraph{Proof of Lemma \ref{lem:conc}: }
Define for $n\geq 1$ the intervals $I_n = [u - 2^{n-1} \eta, u + 2^{n-1} \eta]$ and let $M$ and $K_0$ be sufficiently
large numbers. We have $[-K_0 , K_0 ]\subset I_{n_0}$ with $n_0 = C_1 \ln(K_0 /\eta) \leq C_2 \ln(N K_0 )$ for some constants $C_1, C_2>0$. Denote by $A$ the event: 
$$A=\Big \{\max_{n\leq n_o} \frac{\mathcal{N}(I_n)}{\sqrt{Np}  |I_n|}\geq M\Big \}\bigcup \Big \{ \max y_i\geq K_0\Big \},$$
where for a given interval $I$, we denote by $|I|$ its length and set $\mathcal{N}(I)=\sharp\{i=1, \ldots,N| \,y_i\in I\}$.
Let then $\mbP_k$ denote the probability w.r.t. $C_k$. One has that $$ P[|X^{(k)}| \geq \delta] \leq \mbE (1_{A^c}  \mbP_k [|X| \geq \delta]) + P(A).$$
The first term can be handled as in \cite{SchleinYau1}. It is proved in \cite{SchleinYau1}, Proposition 4.5, under assumption $({\bf A_1})$, and if $\eta \geq 1/N $, and for sufficiently large $M$ and $ K_0$, that there is a positive constant $c = c(M, K_0 )$ such that for any $\delta > 0$
$$\mbE \left (1_{A^c} \mbP_k(|X^{(k)}|\geq \delta)\right )\leq 4e^{-c\min \{\delta \sqrt{N\eta}/\sqrt{p/N}, \delta^2 N^2\eta/p\}}.$$
The above estimate requires again the full assumption $(H_3)$.\\
We now turn to the estimate of $\mbP(A)$. This is the only part which has to be modified to consider sample covariance matrices. For this result, we need that $(\ln N )^2 /N \leq \eta \leq 1$. Then for
sufficiently large $M$ and $K_0$ there are positive constants $c,$ $C$ such that for all $N \geq 2$,
\be \label{estP(A)}\mbP (A) \leq Ce^{-c\sqrt{M \eta\sqrt{Np}}}.\ee
To prove (\ref{estP(A)}), we first recall Lemma 7.3 of \cite{SchleinYau0}. Let $Y$ be a $N\times p, p\geq N $ random matrix with i.i.d. centered entries with variance $1$. If the entries of $Y$ also satisfy assumption 
$({\bf A_1})$, there exists a positive constant $c$ such that for $C>0$ large enough 
$$\mbP \left ( y_{max}(YY^*/p)\geq C\right) \leq e^{-cCp}.$$ 
Thus we only have to consider 
$$\mbP\left ( \Big \{\max_{n\leq n_o} \frac{\mathcal{N}(I_n)}{\sqrt{Np}  |I_n|}\geq M\Big \}\right).$$
To this aim we consider an interval $I$ of length $|I|=\alpha$ for some $\alpha\geq \eta$ and call $u$ its midpoint.
In the following we set $z=u+i\alpha$ and assume that $|z|\geq \varepsilon$ for some $\varepsilon >0$ small. 
Let $0<\vartheta_o <1/8$ be given.
For each $k = 1, 2, \ldots , N$ , we define the events
$$B^{(k)}:=\sum_{m, y_{m}\in I}\xi_m^{(k)}\leq \vartheta_o (\mathcal{N}(I)-1),$$ and set $B=\cup_k B^{(k)}.$ 
As the eigenvalues of $H_N$ and $H_N^{(k)}$ are interlaced, at least $\mathcal{N}(I)-1$ eigenvalues of $H_N^{(k)}$ lie in $I$.
It is proved in Lemma 4.7 of \cite{SchleinYau1} that $$\exists c>0, \,\mbP_k( B^{(k)})\leq e^{-c\sqrt{\mathcal{N}(I)-1}}.$$ 
From this, we deduce that there exists $C>0$ such that if $M$ is large enough,
$$\mbP\left ( B \cap \{\mathcal{N}(I)\geq M\alpha \sqrt{Np}\}\right)\leq Ce^{-c \sqrt{M\alpha \sqrt{Np}}}.$$
Then on $B^c$ we have for some constant $C_2>0$
\begin{eqnarray}&\mathcal{N}(I)&\leq C_2N \alpha \Im \left (m_N(u+i\alpha)\right)\crcr
&&=C_2N \alpha \Im \left (\frac{p/N-1}{z}-\frac{1}{N}\sum_{k=1}^p \frac{1}{z+zC_K^* \mathbf{R}_N^{(k)}C_k} \right)\crcr
&&\leq C_2N \alpha \Im\left (-\frac{1}{N}\sum_{k=1}^p \frac{1}{z+zC_K^* \mathbf{R}_N^{(k)}C_k} \right)\label{impos}\\
&& \leq \frac{ 4C_2Np \alpha^2}{\varepsilon \vartheta_o \mathcal{N}(I)}.\label{last line}
\end{eqnarray}
In (\ref{impos}) we used that $\Im (1/z)<0$ and $p/N \geq 1$. To derive (\ref{last line}), we used the fact that
$$\Im\left (-\frac{1}{z+zC_K^* \mathbf{R}_N^{(k)}C_k} \right)\leq \frac{1}{|z| \Big |1+C_K^* \mathbf{R}_N^{(k)}C_k\Big |}
\leq \frac{1}{|z|\Im (C_K^* \mathbf{R}_N^{(k)}C_k)},$$ with $\frac{1}{|z|}\leq \varepsilon^{-1}$ 
and we estimated from below 
$$\sum_{m=1}^N \frac{\xi_m^{k}}{|y_m^{(k)}-z|^2}\geq \vartheta_o \mathcal{N}(I)/(4\alpha^2 ).$$
This now implies that $\mathcal{N}(I)\leq \alpha \sqrt{4C_2Np/\vartheta_o \varepsilon }$ on $B^c$. This yields the desired result : one simply chooses $M$ large enough.
$\square$


\section{Appendix : Asymptotics of Bessel functions \label{app: B}}
We use two types of asymptotics for Bessel functions. \\
The first one deals with is well-known asymptotics of Bessel functions with bounded order and large argument (see \cite{Olver} e.g). They are used in Section \ref{sec: asanal}.

Assume that $\nu$ is bounded. Then, for large $|z|, z\in \mathbb{C}$
\begin{eqnarray} \label{asyBessel}
&& I_{\nu}(z)=\frac{1}{\sqrt{2\pi z}} ( e^z +e^{-z+(\nu+1/2) i\pi } ) (1+O(1/z)), -\frac{\pi}{2} <\text{ Arg }z<\frac{3\pi}{2};\crcr
&& K_{\nu}(z)=\frac{\sqrt \pi}{\sqrt{2 z}}e^{-z}(1+O(1/z)), |\text{Arg z}|\leq \frac{3\pi}{2}.
\end{eqnarray}

We also make use in Section \ref{sec: gennu} of asymptotics of Bessel functions with large order and large argument.
Abramowitz and Stegun \cite{AbSteg} (p 378) give the following uniform asymptotic expansion of modified Bessel functions of large order :
\begin{eqnarray}\label{asyBessel2}
&&I_{\nu}(\nu z)=\frac{e^{\nu \varphi}}{\sqrt{2\pi \nu}(1+z^2)^{1/4}}\left (1+\sum_{k=1}^{\infty}\frac{u_k(t)}{\nu^k} \right), \crcr
&&K_{\nu}(\nu z)=\frac{\sqrt{\pi}e^{-\nu \varphi}}{\sqrt{2 \nu}(1+z^2)^{1/4}}\left (1+\sum_{k=1}^{\infty}\frac{(-1)^ku_k(t)}{\nu^k} \right),\crcr
&& t=\frac{1}{\sqrt{1+z^2}}, \: \varphi = \sqrt{1+z^2}+\ln \left (\frac{z}{1+\sqrt{1+z^2}}\right), |\arg (z)|\leq \frac{\pi}{2}-\epsilon,\crcr
&&
\end{eqnarray}
where $u_k(t)=t^k v_k(t)$ for some polynomial $v_k.$

\end{document}